\documentclass[twoside,11pt]{article}

\usepackage{jmlr2e}

\usepackage{color}
\usepackage{latexsym}
\usepackage{amsmath,amssymb,amsfonts}
\usepackage{algorithm,algorithmic}
\usepackage{graphicx,subfigure}
\usepackage[breaklinks]{hyperref}
\usepackage{enumitem}

\newcommand\myparagraph[1]{\noindent {\it #1.}}
\newcommand\inner[1]{\langle #1 \rangle}
\newcommand{\normL}{| \! | \! |}
\newcommand{\normR}{| \! | \! |}
\newcommand\mNorm[1]{\normL #1 \normR}

\def\Expect{\mathbb{E}}
\def\reals{\mathbb{R}}
\def\Parti{D}
\def\dParti{\Parti'}
\def\ddParti{\Parti''}
\def\dddParti{\Parti'''}
\def\xs{X_r^{(i)}}
\def\xS{X_S^{(i)}}
\def\xt{X_t^{(i)}}
\def\xsc{X_{V\setminus r}^{(i)}}
\def\E{\mathbb{E}}
\def\grad{\nabla}
\def\ConstMu{\kappa_{m}}
\def\ConstVar{\kappa_{v}}
\def\ConstHess{\kappa_{h}}

\def\ConstGLM{\kappa_4}
\def\Cmin{\rho_{\min}}
\def\Dmax{\rho_{\max}}

\def\DConstOne{\kappa_1(n,p)}
\def\DConstTwo{\kappa_2(n,p)}
\def\DConstThree{\kappa_3(n,p)}

\DeclareMathOperator*{\argmin}{arg\,min}
\DeclareMathOperator*{\minimize}{minimize}
\DeclareMathOperator*{\maximize}{maximize}

\ShortHeadings{Graphical Models via Univariate Exponential Family Distributions}{Yang, Ravikumar, Allen and Liu}
\firstpageno{1}

\begin{document}

\title{Graphical Models via Univariate\\Exponential Family Distributions}

\author{\name Eunho Yang \email eunhyang@us.ibm.com \\
	\addr IBM T.J. Watson Research Center\\
	Yorktown Heights, NY 10598, USA
	\AND
	\name Pradeep Ravikumar \email pradeepr@cs.utexas.edu \\
	\addr Department of Computer Science\\
	University of Texas, Austin\\
	Austin, TX 78712, USA
	\AND
	\name Genevera I. Allen \email gallen@rice.edu \\
	\addr Department of Statistics \\
	Rice University\\
	Houston, TX 77005, USA
	\AND
	\name Zhandong Liu \email zhandonl@bcm.edu \\
	\addr Department of Pediatrics-Neurology \\
	Baylor College of Medicine\\
	Houston, TX 77030, USA
}

\editor{Tommi Jaakkola}

\maketitle

\begin{abstract}%   <- trailing '%' for backward compatibility of .sty file
	Undirected graphical models, or Markov networks, are a popular class of statistical models, used in a wide variety of applications. Popular instances of this class include Gaussian graphical models and Ising models. In many settings, however, it might not be clear which subclass of graphical models to use, particularly for non-Gaussian and non-categorical data. In this paper, we consider a general sub-class of graphical models where the node-wise conditional distributions arise from exponential families. This allows us to derive \emph{multivariate} graphical model distributions from \emph{univariate} exponential family distributions, such as the Poisson, negative binomial, and exponential distributions. Our key contributions include a class of M-estimators to fit these graphical model distributions; and rigorous statistical analysis showing that these M-estimators recover the true graphical model structure exactly, with high probability.  We provide examples of genomic and proteomic networks learned via instances of our class of graphical models derived from Poisson and exponential distributions.
\end{abstract}


%%%%%%%%%%%%%%%%%%%%%%%%%%%%%%%%%%%%%%%%%%%%%%%%%%%%%%%
\section{Introduction}

Undirected graphical models, also known as Markov random fields, are
an important class of statistical models that have been extensively
used in a wide variety of domains, including statistical physics,
natural language processing, image analysis, and medicine. The key
idea in this class of models is to represent the joint distribution as
a product of \emph{clique-wise compatibility functions}.  Given an
underlying graph, each of these compatibility functions depends only
on a subset of variables within any clique of the underlying
graph. Popular instances of such graphical models include Ising and
Potts models (see references in \cite{wainwright_jordan_2008} for a
varied set of applications in computer vision, text analytics, and
other areas with discrete variables), as well as Gaussian Markov
Random Fields (GMRFs), which are popular in many scientific settings
for modeling real-valued data. A key modeling question that
arises, however, is: how do we pick the clique-wise compatibility functions, or
alternatively, how do we pick the form or sub-class of the graphical
model distribution (e.g. Ising or Gaussian MRF)? For the case of
discrete random variables, Ising and Potts models are
popular choices; but these are not best suited for count-valued variables, where the values
taken by any variable could range over the entire set of positive
integers. Similarly, in the case of continuous variables, Gaussian
Markov Random Fields (GMRFs) are a popular choice; but the distributional
assumptions imposed by GMRFs are quite stringent. The marginal distribution of any variable
would have to be Gaussian for instance, which might not hold in instances when the random
variables characterizing the data are skewed \citep{liu_nonpara_2009}. More
generally, Gaussian random variables have thin tails, which might not capture fat-tailed
events and variables. For instance, in the finance domain, the lack of
modeling of fat-tailed events and probabilities has been suggested as
one of the causes of the 2008 financial crisis \citep{acemoglu2009crisis}.

To address this modeling question, some have recently proposed non-parametric extensions of graphical models. Some, such as the non-paranormal \citep{liu_nonpara_2009, lafferty_nonpara_2012}
and copula-based methods \citep{ dobra_2011, liu_copula_2012}, use or learn transforms that Gaussianize the data, and then fit Gaussian MRFs to estimate network structure.  Others, use non-parametric
approximations, such as rank-based estimators, to the correlation matrix, and then fit a Gaussian MRF \citep{xue2012regularized, liu_skeptic_2012}.  More broadly, there could be non-parametric methods that either learn the sufficient statistics functions, or learn transformations of the variables, and then fit standard MRFs over the transformed variables. However, the sample complexity of such classes of non-parametric methods is typically inferior to those that learn parametric models. Alternatively, and specifically for the case of multivariate count data, \citet{lauritzen_1996, bishop_discrete_2007} have suggested combinatorial approaches to fitting graphical models, mostly in the context of
contingency tables. These approaches, however, are computationally intractable for even moderate numbers of variables.

%% One approach to address this model selection question is the recently
%% proposed non-paranormal~\citep{liu_nonpara_2009}, which first
%% Gaussianizes the data, and then fits a GMRF over the transformed
%% variables; as well as related copulas based methods [ref]. More
%% broadly, there could be non-parametric methods that either learn the
%% sufficient statistics functions, or learn transformations of the
%% variables, and then fit standard MRFs over the transformed variables;
%% however the sample complexity of such non-parametric methods is
%% typically inferior to those that learn parametric
%% models. Alternatively, specifically for the case of multivariate count
%% data, \citet{lauritzen_1996, bishop_discrete_2007,
%% wainwright_jordan_2008, esl_2nd} have suggested combinatorial
%% approaches to fitting graphical models, mostly in the context of
%% contingency tables. These approaches, however, are computationally
%% intractable for even moderate numbers of variables.

Interestingly, for the case of \emph{univariate} data, we have a good
understanding of appropriate statistical models to use. In particular,
a count-valued random variable can be modeled using a Poisson
distribution; call-times, time spent on websites, diffusion processes,
and life-cycles can be modeled with an exponential distribution; other
skewed variables can be modeled with gamma or chi-squared
distributions. Here, we ask if we can extend this modeling toolkit
from univariate distributions to multivariate graphical model
distributions?  Interestingly, recent state of the art methods for
learning Ising and Gaussian MRFs~\citep{Meinshausen06,RWL10,JRVS11}
suggest a natural procedure deriving such multivariate graphical models from
univariate distributions. The key idea in these recent methods is to
learn the MRF graph structure by estimating node-neighborhoods, or by
fitting node-conditional distributions of
each node conditioned on the rest of the nodes. Indeed, these
node-wise fitting methods have been shown to have strong computational as
well as statistical guarantees. Here, we consider the general class of
models obtained by the following construction: suppose the
node-conditional distributions of each node conditioned on the rest of
the nodes follows a univariate exponential family. By the
Hammersley-Clifford Theorem~\citep{lauritzen_1996}, and some algebra
as derived in \cite{besag_1974}, these node-conditional distributions
entail a global multivariate distribution that (a) factors according to
cliques defined by the graph obtained from the node-neighborhoods, and
(b) has a particular set of compatibility functions specified by the
univariate exponential family. The resulting class of MRFs, which we
call exponential family MRFs, broadens the class of models available
off the shelf, from the standard Ising, indicator-discrete, and
Gaussian MRFs.

Thus the class of exponential family MRFs provides a principled approach to model multivariate distributions and network structures among a large number of variables; by providing a natural
way to ``extend'' univariate exponential families of distributions to the multivariate case, in many cases where multivariate extensions did not exist in an analytical or computationally tractable form.   Potential applications for these exponential family graphical models abound.  Networks of
call-times, time spent on websites, diffusion processes, and
life-cycles can be modeled with exponential graphical models; other
skewed multivariate data can be modeled with gamma or chi-squared
graphical models; while multivariate count data such as from website
visits, user-ratings, crime and disease incident reports, and
bibliometrics could be modeled via Poisson graphical models. A key
motivating application for our research is multivariate count data
from next-generation genomic 
sequencing technologies~\citep{mortazavi2008mapping}. This technology produces read counts of the number of short RNA fragments that have been mapped back to a
particular gene; and measures gene expression with less technical variation than, and is thus rapidly replacing, microarrays~\citep{marioni_2008}. Univariate count data is typically modeled using Poisson or negative binomial distributions~\citep{li_rna_2011}. As Gaussian graphical models have been traditionally used to understand genomic relationships and estimate regulatory pathways from
microarray data, Poisson and negative-binomial graphical models could thus be used to analyze this next-generation sequencing data. Furthermore, there is a proliferation of new technologies to measure
high-throughput genomic variation in which the data is not even approximately Gaussian (single nucleotide polymorphisms, copy number, methylation, and micro-RNA and gene expression via next-generation
sequencing).  For this data, a more general class of high-dimensional graphical models could thus lead to important breakthroughs in understanding genomic relationships and disease networks.

The construction of the class of exponential family graphical models also
suggests a natural method for fitting such models: node-wise
neighborhood estimation via sparsity constrained node-conditional
likelihood maximization. A main contribution of this paper is to
provide a sparsistency analysis (or analysis of variable selection
consistency) for the recovery of the underlying
graph structure of this broad class of MRFs. We note that the presence
of non-linearities arising from the generalized linear models (GLM) posed subtle technical issues
not present in the linear case~\citep{Meinshausen06}. Indeed, for the
specific cases of logistic, and multinomial respectively,
\cite{RWL10,JRVS11} derive such a sparsistency analysis via fairly
extensive arguments, but which were tuned to the specific cases; for
instance they used the fact that the variables were bounded, and the
specific structure of the corresponding GLMs. Here we generalize their
analysis to general GLMs, which required a subtler analysis as well
as a slightly modified M-estimator. We note that this analysis might be of
independent interest even outside the context of modeling and
recovering graphical models. In recent years, there has been a trend
towards unified statistical analyses that provide statistical
guarantees for broad classes of models via general theorems
\citep{negahban2010}. Our result is in this vein and provides
structure recovery for the class of sparsity constrained generalized
linear models. We hope that the techniques we introduce might be of
use to address the outstanding question of sparsity constrained
M-estimation in its full generality.

There has been related work on the simple idea above of constructing
joint distributions via specifying node-conditional
distributions. \citet{VV05,Varin11} propose the class of composite
likelihood models where the joint distribution is a function of the
conditional distributions of subsets of nodes conditioned on other
subsets. \cite{besag_1974} discuss such joint distribution
constructions in the context of node-conditional distributions
belonging to exponential families, but for special cases of joint
distributions such as pairwise models. In this paper, we consider the
general case of higher-order graphical models for the joint
distributions, and univariate exponential families for the
node-conditional distributions.  Moreover, a key contribution of the
paper is that we provide tractable $M$-estimators with corresponding
high-dimensional statistical guarantees and analysis for learning
this class of graphical models even under high-dimensional statistical
regimes.

Additionally, we note that a preliminary abridged version of this paper
appeared at \citep{YRAL12}.  In this manuscript, we provide a more in
depth theoretical analysis along with several novel developments.
Particularly, we provide a novel analytic framework on the sparsistency of our
$M$-estimators that provide tighter finite-sample bounds, simpler
proofs, and less restrictive assumptions than that of \citep{YRAL12};
these innovations are discussed further in Section 3.2. 
Further, we highlight and study several instances of our framework,
relating our work to the existing literature on Gaussian MRFs and
Ising models, as well as introducing two novel instances, the Poisson
MRF and Exponential MRF.  For each of these cases, we provide specific
corollaries on conditions necessary for sparsistent recovery of the
underlying graph structure.  Finally, we also provide a greatly
expanded experimental analysis of our class of MRFs and their
$M$-estimators compared to that of \citep{YRAL12}.  Focusing on two
novel instances of our model, the Poisson and Exponential MRF, we
study the theoretical rates, graph structural recovery, and robustness
of our estimators through simulated examples.  Further, we provide an
additional case study on protein signaling networks using the
Exponential MRF in Section 4.2.2.

%%%%%%%%%%%%%%%%%%%%%%%%%%%%%%%%%%%%%%%%%%%%%%%%%%%%%%%
\section{Exponential Family Graphical Models}\label{gmodel}

Suppose $X = (X_1,\hdots,X_p)$ is a random vector, with each variable $X_i$ taking values in a set $\mathcal{X}$. Let $G = (V, E)$ be an undirected graph over the set of nodes $V:=\{1,\hdots,p\}$ corresponding to the $p$ variables $\{X_r\}_{r=1}^{p}$. The graphical model over $X$ corresponding to $G$ is a set of distributions that satisfy \emph{Markov independence assumptions} with respect to the graph $G$~\citep{lauritzen_1996}. By the Hammersley-Clifford theorem~\citep{Clifford90}, any such distribution that is strictly positive over its domain also factors according to the graph in the following way. Let $\mathcal{C}$ be a set of cliques (fully-connected subgraphs) of the graph $G$, and let $\{ \phi_{c}(X_c) \}_{c \in \mathcal{C}}$ be a set of clique-wise sufficient statistics. With this notation, any strictly positive distribution of $X$ within the graphical model family represented by the graph $G$ takes the form:
\begin{eqnarray}
\label{EqnGenMRFPot}
P(X) & \propto & \exp \biggr \{\sum_{c \in \mathcal{C}}
\theta_c \phi_c(X_c) \biggr \},
\end{eqnarray}
where $\{\theta_c\}$ are weights over the sufficient statistics. An important special case is a \emph{pairwise graphical model}, where the set of cliques $\mathcal{C}$ consists of the set of nodes $V$ and the set of edges $E$, so that
\begin{eqnarray}
\label{EqnGenMRFPairwise}
P(X) & \propto & \exp \biggr \{\sum_{r \in V} \theta_r \phi_{r}(X_r) +
\sum_{(r,t) \in E} \theta_{rt} \phi_{rt}(X_r,X_t) \biggr \}.
\end{eqnarray}

As previously discussed, an important question is how to select the
form of the graphical model distribution, which under the above
parametrization in \eqref{EqnGenMRFPot}, translates to the question of
selecting the class of sufficient statistics, $\phi$. As discussed in
the introduction, it is of particular interest to derive such a
graphical model distribution as a \emph{multivariate} extension of
specified \emph{univariate} parametric distributions such as negative
binomial, Poisson, and others.  We next outline a subclass of graphical
models that answer these questions via the simple construction: set
the node-conditional distributions of each node conditioned on the
rest of the nodes as following a univariate exponential family, and then
derive the joint distribution that is consistent with these
node-conditional distributions. Then, in Section 3, we will study how
to learn the underlying graph structure, or the edge set $E$, for this
general class of ``exponential family'' graphical models. We provide a
natural sparsity-encouraging $M$-estimator, and sufficient conditions
under which the $M$-estimator recovers the graph structure with high
probability.

\subsection{The Form of Exponential Family Graphical Models}
A popular class of univariate distributions is the exponential family,
whose distribution for a random variable $Z$ is given by
\begin{align}\label{EqnExpFamDist}
	P(Z) = \exp \Big\{\theta \, B(Z) + C(Z) - D(\theta)\Big\},
\end{align}
with sufficient statistics $B(Z)$, base measure $C(Z)$, and
log-normalization constant $D(\theta)$. Such exponential family
distributions include a wide variety of commonly used distributions
such as Gaussian, Bernoulli, multinomial, Poisson, exponential, gamma,
chi-squared, beta, and many others; any of which can be instantiated
with particular choices of the functions $B(\cdot)$, and
$C(\cdot)$. Such exponential family distributions are thus used to
model a wide variety of data types including skewed continuous data
and count data. Here, we ask if we can leverage this ability to model
univariate data to also model the multivariate case. Let $X =
(X_1,X_2,\hdots,X_p)$ be a $p$-dimensional random vector; and let $G =
(V, E)$ be an undirected graph over $p$ nodes corresponding to the $p$
variables. Could we then derive a graphical model distribution over
$X$ with underlying graph $G$, from a particular choice of
\emph{univariate} exponential family
distribution~\eqref{EqnExpFamDist} above?

Consider the following construction. Set the distribution of $X_r$ given the rest of nodes $X_{V \backslash r}$ to be given by the above univariate exponential family distribution~\eqref{EqnExpFamDist}, and where the canonical exponential family parameter $\theta$ is set to a linear combination of $k$-th order products of univariate functions $\{B(X_t)\}_{t \in N(r)}$, where $N(r)$ is the set of neighbors of node $r$ according to graph $G$. This gives the following conditional distribution:
\begin{align}\label{EqnCondDistK}
	 P(X_r | X_{V \backslash r}) &= \exp
		\bigg\{B(X_r) \Big( \theta_{r} + \sum_{t \in N(r)}
		\theta_{r t} \, B(X_{t})+  \sum_{t_2,t_3 \in N(r)}
		\theta_{r \,t_2 t_3} \, B(X_{t_{2}}) B(X_{t_{3}}) \nonumber\\
		& + \hdots + \sum_{t_{2},\hdots,t_{k} \in N(r)} \theta_{r \, t_2\hdots t_{k}} \,  \prod_{j=2}^{k}B(X_{t_j}) \Big)
	+ C(X_r) - \bar{D}(X_{V \backslash r})\bigg\},
\end{align}
where $C(X_r)$ is specified by the exponential family, and $\bar{D}(X_{V \backslash r})$ is the log-normalization constant. Notice that we use the notation $\bar{D}(\cdot)$ in case when we express the log-partition function in terms of random variables. That is, $\bar{D}(X_{V \backslash r}) := D\big(\theta(X_{V \backslash r}) \big)$ where $\theta(X_{V \backslash r})$ is the canonical parameter $\theta$ derived from $X_{V \backslash r}$.

By the Hammersley-Clifford theorem, and some elementary calculation, this conditional distribution can be shown to specify the following unique joint distribution $P(X_1,\hdots,X_p)$:
\begin{proposition}\label{PropCondJoint}
Suppose $X = (X_1,X_2,\hdots,X_p)$ is a $p$-dimensional random vector, and its node-conditional distributions are specified by \eqref{EqnCondDistK} given an undirected graph $G$. Then its joint distribution
$P(X_1,\hdots,X_p)$ belongs to the graphical model represented by $G$, and is given by:
\vspace{-0.3cm}
\begin{align}\label{EqnJointDistK}
	P(X) &= \exp\Biggr\{\sum_{r \in V} \theta_{r} B(X_r) + \sum_{r \in V}\sum_{t \in N(r)} \theta_{rt} \, B(X_r) B(X_{t}) + \hdots \nonumber\\
	& + \sum_{r \in V}\sum_{t_{2},\hdots,t_{k} \in N(r)} \theta_{r\hdots t_{k}} \,  B(X_r) \prod_{j=2}^{k}B(X_{t_j}) + \sum_{r \in V} C(X_r) - A(\theta)\Biggr\},
\end{align}
where $A(\theta)$ is the log-normalization constant.
\end{proposition}

Note that the function $D(\cdot)$ (and hence $\bar{D}(\cdot)$) in \eqref{EqnCondDistK} is the log-partition function of the node-conditional distribution, while the function $A(\cdot)$ in \eqref{EqnJointDistK} in turn is the log-partition function of the joint distribution. Proposition~\ref{PropCondJoint}, thus, provides an answer to our earlier
question on selecting the form of a graphical model distribution
given a univariate exponential family distribution.  When the
node-conditional distributions follow a univariate exponential
family as in \eqref{EqnCondDistK}, there exists a unique graphical model
distribution as specified by \eqref{EqnJointDistK}. One question that
remains, however, is whether the above construction, beginning with
\eqref{EqnCondDistK}, is the most general
possible.  In particular, note that the canonical parameter of the
node-conditional distribution in \eqref{EqnCondDistK} is a {\em tensor
factorization} of the univariate sufficient statistic, which seems a
bit stringent. Interestingly, by extending the argument from
\citep{besag_1974}, which considers the special pairwise case, and
the Hammersley-Clifford Theorem, we can show that indeed
\eqref{EqnCondDistK} and \eqref{EqnJointDistK} have the most general
form.

\begin{theorem}\label{PropGeneralCondForm}
Suppose $X = (X_1,X_2,\hdots,X_p)$ is a $p$-dimensional random vector, and its node-conditional distributions are specified by an exponential family,
\begin{align}\label{EqnCondDistTmp}
	P(X_r | X_{V \backslash r}) = \exp \Big\{ E(X_{V \backslash r})\, B(X_r) + C(X_r) - \bar{D}(X_{V \backslash r}) \Big\},
\end{align}
where the function $E(X_{V \backslash r})$, the canonical parameter of exponential family, depends on the rest of all random variables except $X_r$ (and hence the log-normalization constant $\bar{D}(X_{V \backslash r})$). Further, suppose the corresponding joint distribution factors according to the graph $G$, with the factors over cliques of size at most $k$. Then, the conditional distribution in \eqref{EqnCondDistTmp} necessarily has the tensor-factorized form in \eqref{EqnCondDistK}, and the corresponding joint distribution has the form in \eqref{EqnJointDistK}.
\end{theorem}

\noindent
Theorem~\ref{PropGeneralCondForm} thus tells us that under the general assumptions that:
%\begin{enumerate}[leftmargin=0.05cm, itemindent=0.85cm,label=\textbf{(\alph*)}]
\begin{description}[leftmargin=0.5cm]
	\item[(a)] the joint distribution is a graphical model that factors according to a graph $G$, and has clique-factors of size at most $k$, and
	\item[(b)] its node-conditional distribution follows an exponential family,
\end{description}
%\end{enumerate}
it \emph{necessarily} follows that the conditional and joint distributions are given by \eqref{EqnCondDistK} and \eqref{EqnJointDistK} respectively.

An important special case is when the joint graphical model distribution has clique factors of size at most two. From Theorem~\ref{PropGeneralCondForm}, the conditional distribution is given by:
\begin{align}\label{EqnCondDistP}
	P(X_r | X_{V \backslash r}) = \exp \bigg\{\theta_{r} \, B(X_r) + \sum_{t \in N(r)} \theta_{r t} \, B(X_r) B(X_{t})
	+ C(X_r) - \bar{D}(X_{V \backslash r}) \bigg\},
\end{align}
while the joint distribution is given as:
\begin{align}\label{EqnJointDistP}
	P(X) = \exp \bigg\{\sum_{r \in V} \theta_{r} B(X_r) + \sum_{(r,t) \in E }\theta_{rt} \, B(X_r) B(X_{t}) + \sum_{r \in V} C(X_r) - A(\theta) \bigg\}.
\end{align}
For many classes of models (e.g. general Ising, discrete CRFs), the log-partition function of the joint distribution, $A(\cdot)$, has no analytical form, and might even be intractable to compute,  while the function $D(\cdot)$ typically is more amenable, and available in analytical form, since it is the log-partition function of a \emph{univariate} exponential family distribution.

When the univariate sufficient statistic function $B(\cdot)$ is a linear function $B(X_r) = X_r$, then the conditional distribution in \eqref{EqnCondDistP} is precisely a generalized linear model~\citep{McCullagh} in canonical form,
\begin{align}\label{EqnCondDistGLM}
	P(X_r | X_{V \backslash r}) = \exp\bigg\{\theta_{r} \, X_r + \sum_{t \in N(r)} \theta_{r t} \, X_r \, X_{t}
	+ C(X_r) - \bar{D}(X_{V \backslash r};\theta)\bigg\},
\end{align}
where the canonical parameter of GLMs becomes $\theta_{r} + \sum_{t \in N(r)} \theta_{r t} \, X_{t}$. At the same time, the joint distribution has the form
\begin{align}\label{EqnJointDistGLM}
	P(X) = \exp\bigg\{\sum_{r \in V} \theta_{r} X_r + \sum_{(r,t) \in E }\theta_{rt} \, X_r \, X_{t} + \sum_{r \in V} C(X_r) - A(\theta)\bigg\},
\end{align}
where the log-partition function $A(\cdot)$ in this case is defined as
\begin{align}\label{EqnA}
    A(\theta) := \log \int_{X} \exp\bigg\{\sum_{r \in V} \theta_{r} X_r + \sum_{(r,t) \in E }\theta_{rt} \, X_r \, X_{t} + \sum_{r \in V} C(X_r) \bigg\} dX.
\end{align}
%where $\nu$ is an underlying measure with respect to which the density in \eqref{EqnJointDistGLM} is taken.

We will now provide some examples of our general class of ``exponential family'' graphical model distributions, focusing on the case in \eqref{EqnJointDistGLM} with linear functions $B(X_r) = X_r$.
For each of these examples, we will also detail the domain, $\Theta := \{\theta: A(\theta) < + \infty\}$, of valid parameters that ensure that the density is normalizable. Indeed, such constraints on valid parameters are typically necessary for the distributions over countable discrete or continuous valued variables.

\paragraph{Gaussian Graphical Models.}
The popular Gaussian graphical model~\citep{Speed86} can be derived as an instance of the construction in Theorem~\ref{PropGeneralCondForm}, with the univariate Gaussian distribution as the exponential family distribution. The univariate Gaussian distribution with known variance $\sigma^2$ is given by,
\[ P(Z) \propto \exp\left\{ \frac{\mu}{\sigma}\frac{Z}{\sigma} - \frac{Z^2}{2\sigma^2} \right\},\]
where $Z \in \reals$, so that it can be seen to be an exponential family distribution of the form~\eqref{EqnExpFamDist}, with sufficient statistic $B(Z) = \frac{Z}{\sigma}$, and base measure $C(Z) = - \frac{Z^2}{2\sigma^2}$. Substituting these in \eqref{EqnJointDistGLM}, we get the distribution,
\begin{align}\label{EqnJointGaussian}
	&P(X;\theta) \propto \exp\Bigg\{ \sum_{r \in V} \frac{1}{\sigma_r}\theta_r X_r + \sum_{(r,t) \in E}\frac{1}{\sigma_r\sigma_t}\theta_{rt} \, X_r \, X_t - \sum_{r \in V}\frac{X_r^2}{2\sigma_r^2} \Bigg\},
\end{align}
which can be seen to be the multivariate Gaussian distribution. Note that the set of parameters $\{\theta_{rt}\}_{(r,t) \in E}$ entails a precision matrix that needs to be positive definite for a valid probability distribution.    

\paragraph{Ising Models}
The Ising model~\citep{wainwright_jordan_2008} in turn can be derived from the construction in Theorem~\ref{PropGeneralCondForm} with the Bernoulli distribution as the univariate exponential family distribution. The Bernoulli distribution is a member of the exponential family of the form~\eqref{EqnExpFamDist}, with sufficient statistic $B(X) = X$, and base measure $C(X) = 0$, and with variables taking values in the set $\mathcal{X} = \{0,1\}$. Substituting these in \eqref{EqnJointDistGLM}, we get the distribution,
\begin{align}\label{EqnJointIsing}
	P(X;\theta) = \exp\Bigg\{\sum_{(r,t) \in E}\theta_{rt} \, X_r \, X_{t} - A(\theta)\Bigg\},
\end{align}
where we have ignored the singleton term, i.e. set $\theta_r = 0$ for simplicity. The form of the multinomial graphical model, an extension of the Ising model, can also be represented by \eqref{EqnJointDistGLM} and has been previously studied in \citet{JRVS11} and others. The Ising model imposes no constraint on its parameters, $\{\theta_{rt}\}$, for normalizability, since there are finitely many configurations of the binary random vector $X$.

\paragraph{Poisson Graphical Models}
Poisson graphical models are an interesting instance with the Poisson distribution as the univariate exponential family distribution. The Poisson distribution is a member of the exponential family of the form~\eqref{EqnExpFamDist}, with sufficient statistic $B(X) = X$ and $C(X) = - \mathrm{log}( X!)$, and with variables taking values in the set $\mathcal{X} = \{0,1,2,...\}$. Substituting these in \eqref{EqnJointDistGLM}, we get the following Poisson graphical model distribution:
\begin{align}\label{EqnPoisson}
	P(X;\theta) = \exp\Bigg\{\sum_{r \in V} \big( \theta_r X_r - \log (
	X_{r}!) \big) + \sum_{(r,t) \in E }\theta_{rt} \, X_r \, X_{t}  -
	A(\theta)\Bigg\}.
\end{align}
For this Poisson family, with some calculation, it can be seen that the
normalizability condition, $A(\theta) < + \infty$, entails $\theta_{rt}
\leq 0 \ \forall \ r,t$. In other words, the Poisson
graphical model can only capture {\it negative} conditional relationships
between variables.

\paragraph{Exponential Graphical Models}\label{SecExpMRF}
Another interesting instance uses the exponential distribution as the univariate exponential family distribution, with sufficient statistic $B(X) = -X$ and $C(X) = 0$, and with variables taking values in $\mathcal{X} = \{0\}\cup \reals^+ $. Such exponential distributions are typically used for data describing inter-arrival times between events, among other applications.  Substituting these in \eqref{EqnJointDistGLM}, we get the following exponential graphical model distribution:
\begin{align}\label{EqnExp}
	P(X;\theta) = \exp\Bigg\{ - \sum_{r \in V}  \theta_r X_r -  \sum_{(r,t)
	\in E }\theta_{rt} \, X_r \, X_{t}  -
	A(\theta)\Bigg\}.
\end{align}
To ensure that the distribution is valid and normalizable, so that
$A(\theta) < + \infty$, we then require that $\theta_r > 0,
\theta_{rt} \geq0 \ \forall \ r,t$. Because of the negative sufficient
statistic, this implies that the exponential graphical model can only
capture {\it negative} conditional relationships between variables.

%%%%%%%%%%%%%%%%%%%%%%%%%%%%%%%%%%%%%%%%%%%%%%%%%%%%%%%
\section{Statistical Guarantees on Learning Graphical Model Structures}
In this section, we study the problem of learning the graph structure of an underlying exponential family MRF, given i.i.d. samples. Specifically, we assume that we are given $n$ samples of random vector $X^{1:n} := \{X^{(i)}\}_{i=1}^{n}$, from a pairwise exponential family MRF,
\begin{align}\label{EqnJointGLM}
	P(X;\theta^*) = \exp\Bigg\{\sum_{r \in V} \theta^*_{r} X_r + \sum_{(r,t) \in E^*}\theta^*_{rt} \, X_r \, X_{t} + \sum_{r} C(X_r) - A(\theta^*)\Bigg\}.
\end{align}
The goal in graphical model structure recovery is to recover the edges $E^*$ of the underlying graph $G = (V,E^*)$. Following \cite{Meinshausen06,RWL10,JRVS11}, we will approach this problem via neighborhood estimation: where we estimate the neighborhood of each node individually, and then stitch these together to form the global graph estimate. Specifically, if we have an estimate $\widehat{N}(r)$ for the true neighborhood $N^*(r)$, then we can estimate the overall graph structure as,
\begin{align}\label{EqnGraphEstimate}
	\widehat{E} = \cup_{r \in V} \cup_{t \in \widehat{N}(r)}\{(r,t)\}.
\end{align}
\paragraph{Remark.} Note that the node-neighborhood estimates $\widehat{N}(r)$ might not be symmetric (i.e. there may be a pair $(r,s) \in V \times V$, with $r \in \widehat{N}(s)$, but $s \not\in \widehat{N}(r)$). The graph-structure estimate in \eqref{EqnGraphEstimate} provides one way to reconcile these neighborhood estimates; see \citet{Meinshausen06} for some other ways to do so (though as they note, these different estimates have asymptotically identical sparsistency guarantees: given exponential convergence in the probability of node-neighborhood recovery to one, the probability that the node-neighborhood estimates are symmetric, and hence that the different ``reconciling'' graph estimates would become identical, also converges to one.)

The problem of graph structure recovery can thus be reduced to the problem of recovering the neighborhoods of all the nodes in the graph. In order to estimate the neighborhood of any node in turn, we consider the sparsity constrained conditional MLE. Note that given the joint distribution in \eqref{EqnJointGLM}, the conditional distribution of $X_r$ given the rest of the nodes is reduced to a GLM and given by,
\begin{align}\label{EqnCondGLM}
	P(X_r|X_{V \backslash r}) = \exp\Bigg\{ X_r \Big(\theta^*_r + \sum_{t \in N^*(r)} \theta^*_{rt} X_t\Big) + C(X_r) - D\Big(\theta^*_r + \sum_{t \in N^*(r)} \theta^*_{rt} X_t\Big)\Bigg\}.
\end{align}
Let $\theta^*(r)$ be a set of parameters related to the node-conditional distribution of node $X_r$, i.e. $\theta^*(r) = (\theta^*_r,\theta^*_{\backslash r}) \in \reals \times \reals^{p-1}$ where $\theta^*_{\backslash r} = \{\theta^*_{rt}\}_{t \in V\backslash r}$ be a zero-padded vector, with entries $\theta^*_{rt}$ for $t \in N^*(r)$ and $\theta^*_{rt} = 0$, for $t \not\in N^*(r)$. In order to infer the neighborhood structure for each node $X_r$, we solve the $\ell_1$ regularized conditional log-likelihood loss:
\begin{align}\label{Eq:Obj}
    \minimize_{\theta(r) \in \Omega} & \ \Big\{\ell\big(\theta(r);X^{1:n}\big) + \lambda_n \|\theta_{\backslash r}\|_1\Big\},
\end{align}
where $\Omega$ is the parameter space in $\reals \times \reals^{p-1}$, and $\ell\big(\theta(r);X^{1:n}\big)$ is the conditional log-likelihood of the distribution~\eqref{EqnCondGLM},
\begin{align*}
    \ell\big(\theta(r);X^{1:n}\big) &:=  -\frac{1}{n}\log \prod_{i=1}^{n} P\big(\xs | \xsc, \theta(r)\big) \nonumber\\
    &= \frac{1}{n} \sum_{i=1}^n \bigg\{- \xs \Big(\theta_r + \inner{\theta_{\backslash r},\xsc}\Big)  + \Parti\Big(\theta_r + \inner{\theta_{\backslash r},\xsc}\Big) \bigg\}.\nonumber
\end{align*}
Note that the parameter space $\Omega$ might be restricted, and strictly smaller than $\reals \times \reals^{p-1}$; for Poisson graphical models, $\theta_{rt} \leq 0$ for all $r,t \in V$ for instance.

Given the solution $\widehat{\theta}(r)$ of the $M$-estimation problem above, we then estimate the node-neighborhood of $r$ as $\widehat{N}(r) = \{t \in V \backslash r \,:\, \widehat{\theta}_{rt} \neq 0\}$. In what follows, when we focus on a fixed node $r \in V$, we will overload notation, and use $\theta \in \reals \times \reals^{p-1}$ as the parameters of the conditional distribution, suppressing dependence on the node $r$.

\subsection{Conditions}
A key quantity in the analysis is the Fisher information matrix, $Q^*_{r} = \grad^2 \ell\big(\theta^*;X^{1:n}\big)$, which is the Hessian of the node-conditional log-likelihood.  In the following, we again will simply use $Q^*$ instead of $Q^*_r$ where the reference node $r$ should be understood implicitly. We also use $S = \{(r,t) : t \in N^*(r)\}$ to denote the true neighborhood of node $r$, and $S^c$ to denote its complement. We use $Q^*_{SS}$ to denote the $d \times d$ sub-matrix of $Q^*$ indexed by $S$ where $d$ is the number of neighborhoods of node $r$ again suppressing dependence on $r$. Our first two conditions, mirroring those in \cite{RWL10}, are as follows.

\begin{enumerate}[leftmargin=0.05cm, itemindent=1cm,label=\textbf{(C\arabic*)}, ref=(C\arabic*),start=1]
    \item (Dependency condition) There exists a constant $\Cmin > 0$ such that $\lambda_{\min}(Q^*_{SS}) \ge \Cmin$ so that the sub-matrix of Fisher information matrix corresponding to true neighborhood has bounded eigenvalues. Moreover, there exists a constant $\Dmax < \infty$ such that $\lambda_{\max}(\frac{1}{n}\sum_{i=1}^n [\xsc (\xsc)^T]) \le \Dmax$.\label{AN_Dep}
\end{enumerate}
These condition can be understood as ensuring that variables do not become overly dependent. We will also need an incoherence or irrepresentable condition on the Fisher information matrix as in \cite{RWL10}.

\begin{enumerate}[leftmargin=0.05cm, itemindent=1cm,label=\textbf{(C\arabic*)}, ref=(C\arabic*),start=2]
    \item (Incoherence condition) There exists a constant $\alpha > 0$, such that $\max_{t \in S^c} \|Q^*_{tS} (Q^*_{SS})^{-1}\|_1 \allowbreak \le 1 - \alpha$.
\end{enumerate}
This condition, standard in high-dimensional analyses, can be understood as ensuring that irrelevant variables do not exert an overly strong effect on the true neighboring variables.

A key technical facet of the linear, logistic, and multinomial models in \cite{Meinshausen06,RWL10,JRVS11}, used heavily in their proofs, was that the random variables $\{X_{r}\}$ there were bounded with high probability. Unfortunately, in the general exponential family distribution in \eqref{EqnCondGLM}, we cannot assume this explicitly. Nonetheless, we show that we can analyze the corresponding regularized M-estimation problems under the following mild conditions on the log-partition functions of the joint and node-conditional distributions.

\begin{enumerate}[leftmargin=0.05cm, itemindent=1cm,label=\textbf{(C\arabic*)}, ref=(C\arabic*),start=3]
    \item (Bounded Moments) For all $r \in V$, the first and second moments are bounded, so that
    \begin{align*}
        \E[X_{r}] \le \ConstMu \quad \textrm{and} \quad \E[X_{r}^2] \le \ConstVar,
    \end{align*}
    for some constants $\ConstMu$, $\ConstVar$. Further, the log-partition function $A(\cdot)$ of the joint distribution~\eqref{EqnJointGLM} satisfies:
    \begin{align*}
        \max_{u : |u| \le 1} \frac{\partial^{2}}{\partial\theta_r^2}A(\theta^* + u e_r) \le \ConstHess,
    \end{align*}
    for some constant $\ConstHess$, and where $e_r \in \reals^{p^2}$ is an indicator vector that is equal to one at the index corresponding to $\theta_r$, and zero everywhere else. Further, it holds that
    \begin{align*}
        \max_{\eta : |\eta| \le 1} \frac{\partial^{2}}{\partial\eta^2}\bar{A}_r(\eta; \theta^*) \le \ConstHess,
    \end{align*}
    where $\bar{A}_r(\eta;\theta^*)$ is a slight variant of \eqref{EqnA}:
    \begin{align}\label{ABar}
        \bar{A}_r(\eta;\theta) := \log \int_{X} \exp \bigg\{ \eta X_r^2 + \sum_{u \in V} \theta_{u} X_{u} + \sum_{(u,t) \in V^2 }\theta_{ut} \, X_{u} \, X_{t} + \sum_{u \in V} C(X_{u}) \bigg\} \, dX,
    \end{align}
    for some scalar variable $\eta$.\label{AN_GLM_MH}
    \item For all $r \in V$, the log-partition function $D(\cdot)$ of the node-wise conditional distribution \eqref{EqnCondGLM} satisfies: there exist functions $\DConstOne$ and $\DConstTwo$ (that depend on the exponential family) such that, for all $\theta \in \Theta$ and $X \in \mathcal{X}$, $|D''(a)| \le \DConstOne$ where $a \in \big[b, b + 4 \DConstTwo \max\{\log n,\log p\}\big]$ for $b := \theta_r + \inner{\theta_{\backslash r},\, X_{V\backslash r}}$.  Additionally, $|D'''(b)| \leq \DConstThree$ for all $\theta \in \Theta$ and $X \in \mathcal{X}$. Note that $\DConstOne$,$\DConstTwo$ and $\DConstThree$ are functions that might be dependent on $n$ and $p$, which affect our main theorem below.\label{AN_GLM_Cond_Smooth}
\end{enumerate} 
Conditions \ref{AN_GLM_MH} and \ref{AN_GLM_Cond_Smooth} are the key technical components enabling us to generalize the analyses in \cite{Meinshausen06,RWL10,JRVS11} to the general exponential family case. It is also important to note that almost all exponential family distributions including all our previous examples can satisfy \ref{AN_GLM_Cond_Smooth} with mild functions $\DConstOne$, $\DConstTwo$ and $\DConstThree$, as we will explicitly show later in this section. 
Comparing to the assumption in \citet{YRAL12} that requires $\|\theta^*\|_2 \leq 1$ for some exponential families, this will be much less restrictive condition on the minimum values of $\theta^{*}$ permitted to achieve
  variable selection consistency.

\subsection{Statement of the Sparsistency Result}
Armed with the conditions above, we can show that the random vector $X$ following a exponential family MRF distribution in \eqref{EqnJointGLM} is suitably well-behaved:
\begin{proposition}\label{PropSubGaussVariableGLM}
Suppose $X$ is a random vector with the distribution specified in \eqref{EqnJointGLM}. Then, for $\forall r \in V$,		
	\begin{align*}
		P\bigg(\frac{1}{n}\sum_{i=1}^n\big(\xs\big)^2 \geq \delta \bigg) \leq \exp\left(-c \, n \, \delta^2 \right)
	\end{align*}
where $\delta \le \min\{2\ConstVar/3,\ConstHess+\ConstVar\}$, and $c$ is a positive constant.
\end{proposition}
We recall the notation that the superscript indicates the sample and the subscript indicates the node; so that $X^{(i)}$ is the i-th sample, while $X_s^{(i)}$ is the $s$-th variable/node of this random vector.

\begin{proposition}\label{PropSubGaussSingleVariableGLM}
Suppose $X$ is a random vector with the distribution specified in \eqref{EqnJointGLM}. Then, for $\forall r \in V$,
	\begin{align*}
		P\Big(|X_r| \geq \delta \log \eta \Big) \leq c \eta^{-\delta}
	\end{align*}
where $\delta$ is any positive real value, and $c$ is a positive constant.
\end{proposition}
These propositions are key to the following sparsistency result for the general family of pairwise exponential family MRFs~\eqref{EqnJointGLM}.

\begin{theorem}\label{Thm:Main}
Consider a pairwise exponential family MRF distribution as specified in \eqref{EqnJointGLM}, with true parameter $\theta^*$ and associated edge set $E^*$ that satisfies Conditions \ref{AN_Dep}-\ref{AN_GLM_Cond_Smooth}.  Suppose that $\min_{(s,t) \in E^*}|\theta^*_{rt}| \geq \frac{10}{\Cmin}\sqrt{d}\lambda_n$, where $d$ is the maximum neighborhood size. Suppose also that the regularization parameter is chosen such that $M_1 \frac{(2-\alpha)}{\alpha}\sqrt{\DConstOne}\sqrt{\frac{\log p}{n}} \leq \lambda_n \leq M_2 \frac{(2-\alpha)}{\alpha}\DConstOne \DConstTwo$ for some constants $M_1, M_2 > 0$. Then, there exist positive constants $L$, $c_1$, $c_2$ and $c_3$ such that if $n \geq L d^{2}\DConstOne (\DConstThree)^2 \log p ( \max\{\log n,\log p\})^{2}$, then with probability at least $1- c_1 (\max\{n,p\})^{-2} - \exp(-c_2 n) - \exp(-c_3 n)$, the following statements hold.
    %\begin{enumerate}[leftmargin=0.05cm, itemindent=0.85cm,label=\textbf{(\alph*)}]
\begin{description}[leftmargin=0.5cm]
	\item[(a)] (Unique Solution) For each node $r \in V$, the solution of the M-estimation problem in \eqref{Eq:Obj} is unique, and
	\item[(b)] (Correct Neighborhood Recovery) The M-estimate also recovers the true neighborhood exactly, so that $\widehat{N}(r) = N^*(r)$.
\end{description}
%\end{enumerate}

%    \begin{enumerate}[(a)]
%        \item (Unique Solution) For each node $r \in V$, the solution of the M-estimation problem in \eqref{Eq:Obj} is unique, and
%        \item (Correct Neighborhood Recovery) The M-estimate also recovers the true neighborhood exactly, so that $\widehat{N}(r) = N^*(r)$.
%    \end{enumerate}
\end{theorem}
Note that if the neighborhood of each node is recovered with high probability, then by a simple union bound, the estimate in \eqref{EqnGraphEstimate},
$\widehat{E} = \cup_{r \in V} \cup_{t \in \widehat{N}(r)}\{(r,t)\}$ is equal to the true edge set $E^*$ with high-probability.

% \paragraph{Remark.}
% It is instructive to compare the sample complexity derived in Theorem~\ref{Thm:Main} to that derived for multinomial graphical models in \cite{JRVS11}, where authors focused on structure recovery for multinomial graphical models, where each variable $X_r$ took values in a set of cardinality $m$, and entailed a sample size of $n \geq L d^{2} (\log p) m^{2}$. In our theorem that applies to general GLMs, we leverage Proposition~\ref{PropSubGaussSingleVariableGLM}, and show that each variable $X_r$ is upper bounded by $\max\{\log n, \log p\}$ with high probability for any GLM graphical model; and match the sample complexity of \cite{JRVS11}, with $\max\{\log n, \log p\}$ instead of $m$.

In the following subsections, we investigate the consequences of Theorem~\ref{Thm:Main} for the sparsistency of specific instances of our general exponential family MRFs.

\subsection{Statistical Guarantees for Gaussian MRFs, Ising Models, Exponential Graphical Models}
In order to apply Theorem~\ref{Thm:Main} to a specific instance of our general exponential family MRFs, we need to specify the terms $\DConstOne$, $\DConstTwo$ and $\DConstThree$ defined in Condition \ref{AN_GLM_Cond_Smooth}. It turns out that we can specify these terms for the Gaussian graphical models, Ising models and Exponential graphical model distributions, discussed in Section~\ref{gmodel}, in a similar manner, since the node-conditional log-partition function $D(\cdot)$ for all these distributions can be upper bounded by some constant independent of $n$ and $p$. In particular, we can set $\DConstOne := \kappa_1$, $\DConstTwo := \infty$ and $\DConstOne := \kappa_3$ where $\kappa_1$ and $\kappa_3$ now become some constants depending on the distributions.

\paragraph{Gaussian MRFs.}
Recall that the node-conditional distribution for Gaussian MRFs follow a univariate Gaussian distribution :
\begin{align*}%\label{EqnCondGaussian}
	&\hspace{-.5cm}P(X_r|X_{V \backslash r}) \propto \exp\bigg\{ X_r \big(\theta_r + \sum_{t \in N(r)} \theta_{rt} X_t\big) - \frac{1}{2}X_r^2 - \frac{1}{2}\big(\theta_r + \sum_{t \in N(r)} \theta_{rt} X_t\big)^2 \bigg\}.
\end{align*}
Note that following \citep{Meinshausen06}, we assume that $\sigma_r^2 =1$ for all $r \in V$. The node-conditional log-partition function $D(\cdot)$ can thus be written as $D(\eta) := -\frac{1}{2}\eta^2$, so that $|D''(\eta)| = 1$ and $D'''(\eta) = 0$. We can thus set $\kappa_1 = 1$ and $\kappa_3 =0$.

\paragraph{Ising Models.}
    For Ising models, node-conditional distribution follows a Bernoulli distribution:
    \begin{align*}%\label{EqnCondIsing}
    	P(X_r|X_{V \backslash r}) = \exp\Bigg\{ X_r \Big(\sum_{t \in N(r)} \theta_{rt} X_t\Big) - \log \bigg(1 + \exp\Big(\sum_{t \in N(r)} \theta_{rt} X_t\Big)\bigg)\Bigg\}.
    \end{align*}
The node-conditional log-partition function  $D(\cdot)$ can thus be written as $D(\eta) := \log \big(1 + \exp (\eta) \big)$, so that for any $\eta$, $|D''(\eta)| = \frac{\exp(\eta)}{(1+\exp(\eta))^2} \leq \frac{1}{4}$ and $|D'''(\eta)| =  \big|\frac{\exp(\eta)(1-\exp(\eta))}{(1+\exp(\eta))^3}\big| < \frac{1}{4}$. Hence, we can set $\kappa_1 = \kappa_3 = 1/4$.

\paragraph{Exponential Graphical Models.}
    Lastly, for exponential graphical models, we have
    \begin{align*}
    	P(X_r|X_{V \backslash r}) = \exp\Bigg\{ -X_r \Big(\theta_r + \sum_{t \in N(r)} \theta_{rt} X_t\Big) + \log\Big(\theta_r + \sum_{t \in N(r)} \theta_{rt} X_t\Big) \Bigg\}.
    \end{align*}
The node-conditional log-partition function  $D(\cdot)$ can thus be written as $D(\eta) := - \log \eta$, with $\eta = \theta_r + \sum_{t \in N(r)} \theta_{rt} X_t$. Recall from Section~\ref{SecExpMRF} that the node parameters are strictly positive $\theta_r > 0$, and the edge-parameters are positive as well, $\theta_{rt} \geq 0$, as are the variables themselves $X_t \ge 0$. Thus, under the additional constraint that $\theta_r > a_0$ where $a_0$ is a constant smaller than $\theta^*_r$, we have that $\eta := \theta_r + \sum_{t \in N(r)} \theta_{rt} X_t \ge a_0$. Consequently, $|D''(\eta)| = \frac{1}{\eta^2} \leq \frac{1}{a_0^2}$ and $|D'''(\eta)| = \big|\frac{2}{\eta^3}\big| \leq \frac{2}{a_0^3}$. We can thus set $\kappa_1 = \frac{1}{a_0^2}$ and $\kappa_3 =\frac{2}{a_0^3}$.

Armed with these derivations, we recover the following result on the sparsistency of Gaussian, Ising and Exponential graphical models, as a corollary of Theorem~\ref{Thm:Main}:
\begin{corollary}\label{Cor:Ising}
Consider a Gaussian MRF \eqref{EqnJointGaussian} or Ising model \eqref{EqnJointIsing} or Exponential graphical model \eqref{EqnExp} distribution with true parameter $\theta^*$, and associated edge set $E^*$, and which satisfies Conditions~\ref{AN_Dep}-\ref{AN_GLM_MH}.  Suppose that $\min_{(s,t) \in E^*}|\theta^*_{rt}| \geq \frac{10}{\Cmin}\sqrt{d}\lambda_n$. Suppose also that the regularization parameter is set so that $M \frac{(2-\alpha)}{\alpha}\sqrt{\kappa_1}\sqrt{\frac{\log p}{n}} \leq \lambda_n$ for some constant $M > 0$. Then, there exist positive constants $L$, $c_1$, $c_2$ and $c_3$ such that if $n \geq L \kappa_1 \kappa_3^2 d^{2} \log p ( \max\{\log n,\log p\})^{2}$, then with probability at least $1- c_1 (\max\{n,p\})^{-2} - \exp(-c_2 n) - \exp(-c_3 n)$, the statements on the uniqueness of the solution and correct neighborhood recovery, in Theorem \ref{Thm:Main} hold.
\end{corollary}

\paragraph{Remarks.}
As noted, our models and theorems are quite
general, extending well beyond the popular Ising and Gaussian graphical
models. The graph structure recovery problem for Gaussian models was
studied in \cite{Meinshausen06} especially for the regime where the
neighborhood sparsity index is \emph{sublinear}, meaning that $d/p
\rightarrow
0$. Besides the sublinear scaling regime, Corollary \ref{Cor:Ising} can be
adapted to entirely different types of scaling, such as the linear
regime where $d/p
\rightarrow \alpha$ for some $\alpha >0$ (see
\cite{Wainwright2006new} for details on adaptations to sublinear scaling
regimes). Moreover, with $\kappa_1$ and $\kappa_3$ as defined above,
Corollary \ref{Cor:Ising} exactly recovers the result in
\cite{RWL10} for the Ising models as a special case.

Also note that Corollary \ref{Cor:Ising} provides tighter finite-sample bounds than the results of \citet{YRAL12}. In particular, a sample size complexity necessary on $\lambda_n$ to achieve sparsistent recovery here is $O(\sqrt{\frac{\mathrm{log}p}{n}})$, which is faster as compared to $O(\sqrt{\frac{\mathrm{log}p}{n^{1-\kappa}}})$ in \citet{YRAL12}.

\subsection{Statistical Guarantees for Poisson Graphical Models}
We now consider the Poisson graphical model. Again, to derive the corresponding corollary of Theorem~\ref{Thm:Main}, we need to specify the terms $\DConstOne$, $\DConstTwo$ and $\DConstThree$ defined in Condition~\ref{AN_GLM_Cond_Smooth}. Recall that the node-conditional distribution of Poisson graphical models has the form:
\begin{align*}
	P(X_r|X_{V \backslash r}) = \exp\bigg\{X_r \Big(\theta_r + \sum_{t \in N(r)} \theta_{rt} X_t\Big) - \log(X_{r}!) - \exp\Big(\theta_r + \sum_{t \in N(r)} \theta_{rt} X_t\Big) \bigg\}.
\end{align*}
The node-conditional log-partition function  $D(\cdot)$ can thus be written as $D(\eta) := \exp \eta$, with $\eta = \theta_r + \sum_{t \in N(r)} \theta_{rt} X_t$. Noting that the variables $\{X_t\}$ range over positive integers, and that feasible parameters $\theta_{rt}$ are negative, we obtain
\begin{align*}
	D''(\eta) = D''\big(\theta_r + \inner{\theta_{\backslash r},\, X_{V\backslash r}} + 4 \DConstTwo \log p'\big) &= \exp\big(\theta_r + \inner{\theta_{\backslash r},\, X_{V\backslash r}} + 4 \DConstTwo \log p'\big)\\
	& \leq \exp\big(\theta_r + 4 \DConstTwo \log p'\big),
\end{align*}
where $p' = \max\{n,p\}$. Suppose that we restrict our attention on the subfamily where $\theta_r \leq a_0$ for some positive constant $a_0$. Then, if we choose $\DConstTwo := 1/(4\log p')$, we then obtain $\theta_r + 4 \DConstTwo \log p' \le a_0 + 1$,
so that setting $\DConstOne := \exp(a_0 + 1)$ would satisfy Condition~\ref{AN_GLM_Cond_Smooth}.
Similarly, we obtain $D'''\big(\theta_r + \inner{\theta_{\backslash r},\, X_{V\backslash r}}\big) = \exp\big(\theta_r + \inner{\theta_{\backslash r},\, X_{V\backslash r}}\big) \leq \exp(a_0 + 1)$, so that we can set $\DConstThree$ to $\exp(a_0 + 1)$.

Armed with these settings, we recover the following corollary for Poisson graphical models:
\begin{corollary}\label{Cor:Poisson}
Consider a Poisson graphical model distribution as specified in \eqref{EqnPoisson}, with true parameters $\theta^*$, and associated edge set $E^*$, that satisfies Conditions \ref{AN_Dep}-\ref{AN_GLM_MH}. Suppose that $\min_{(s,t) \in E^*}|\theta^*_{rt}| \geq \frac{10}{\Cmin}\sqrt{d}\lambda_n$. Suppose also that the regularization parameter is chosen such that
$M_1 \frac{(2-\alpha)}{\alpha}\sqrt{\kappa_1}\sqrt{\frac{\log p}{n}} \leq \lambda_n \leq M_2 \kappa_1 \frac{(2-\alpha)}{\alpha}\frac{1}{\max\{\log n,\log p\}}$ for some constants $M_1, M_2 > 0$. Then, there exist positive constants $L$, $c_1$, $c_2$ and $c_3$ such that if $n \geq L d^{2}\kappa_1 \kappa_3^2 \log p ( \max\{\log n,\log p\})^{2}$, then with probability at least $1- c_1 (\max\{n,p\})^{-2} - \exp(-c_2 n) - \exp(-c_3 n)$, the statements on the uniqueness of the solution and correct neighborhood recovery, in Theorem \ref{Thm:Main} hold.
\end{corollary}

%%%%%%%%%%%%%%%%%%%%%%%%%%%%%%%%%%%%%%%%%%%%%%%%%%%%%%%
\section{Experiments}

We evaluate our M-estimators for exponential family graphical models,
specifically for the Poisson and exponential distributions,
through simulations and real data examples.  Neighborhood selection
was performed for each M-estimator with an $\ell_{1}$ penalty to
induce sparsity and non-negativity or non-positivity constraints to
enforce appropriate restrictions on the parameters.  Optimization
algorithms were implemented using projected gradient descent
\citep{daubechies2008accelerated,beck_teb_2010},
which since the objectives are convex, is guaranteed to converge to
the global optimum.  Further details on the optimization problems used
for our M-estimators are given in the
Appendix~\ref{appendix:optimization}.

\begin{figure}[h!tb]
    \centering
    \subfigure[Exponential Grid Structure]{\label{Fig:4nn_exponential}
        \centering
        \includegraphics[width=0.4\textwidth]{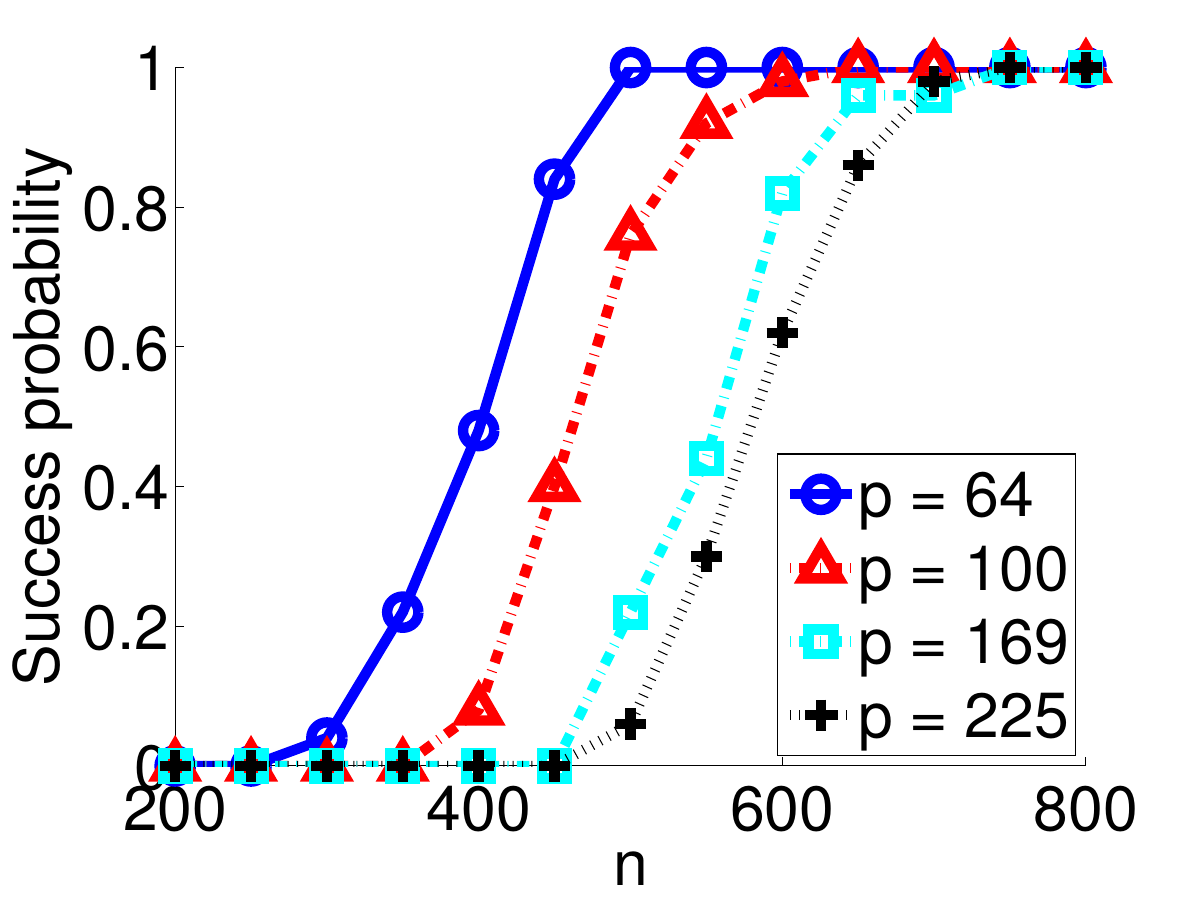}
        \includegraphics[width=0.4\textwidth]{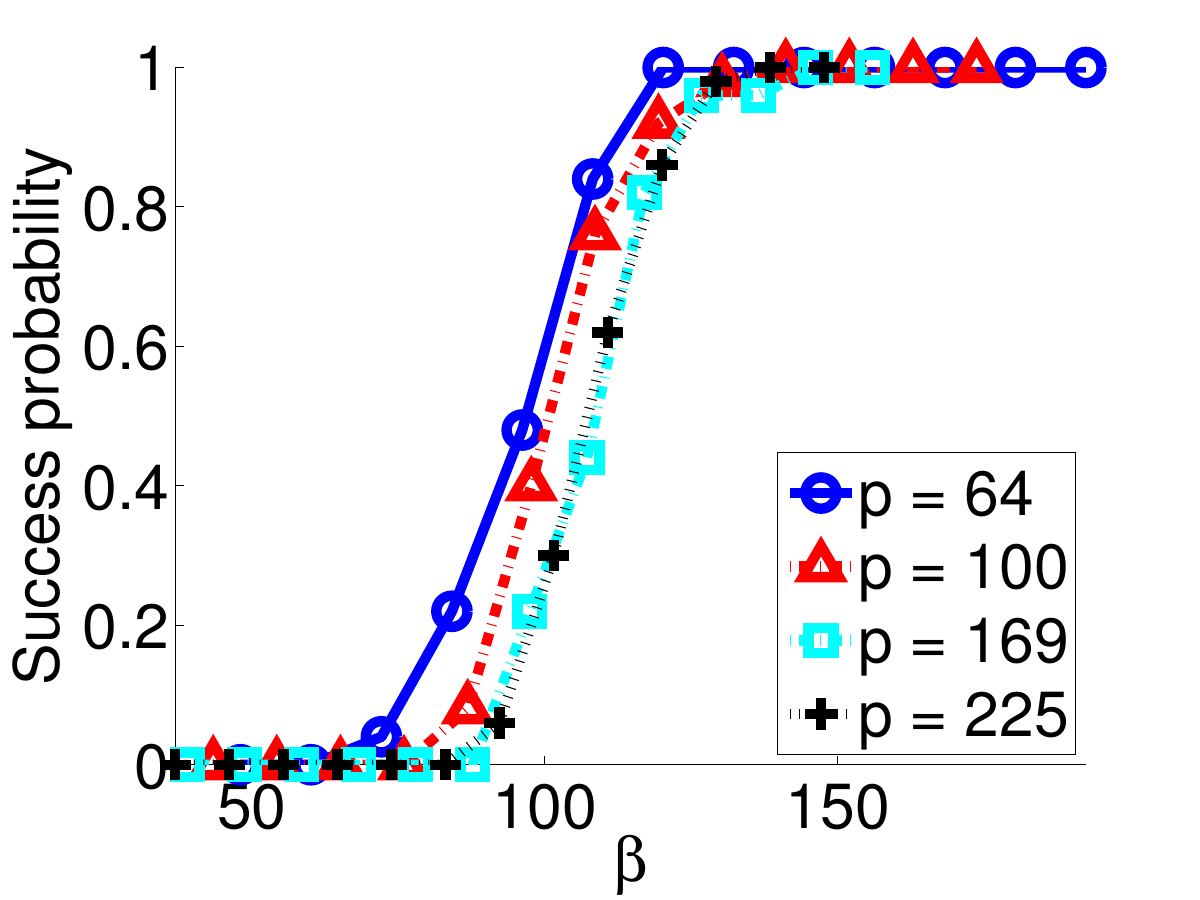}
    }
    \subfigure[Poisson Grid Structure]{\label{Fig:4nn_poisson}
        \centering
        \includegraphics[width=0.4\textwidth]{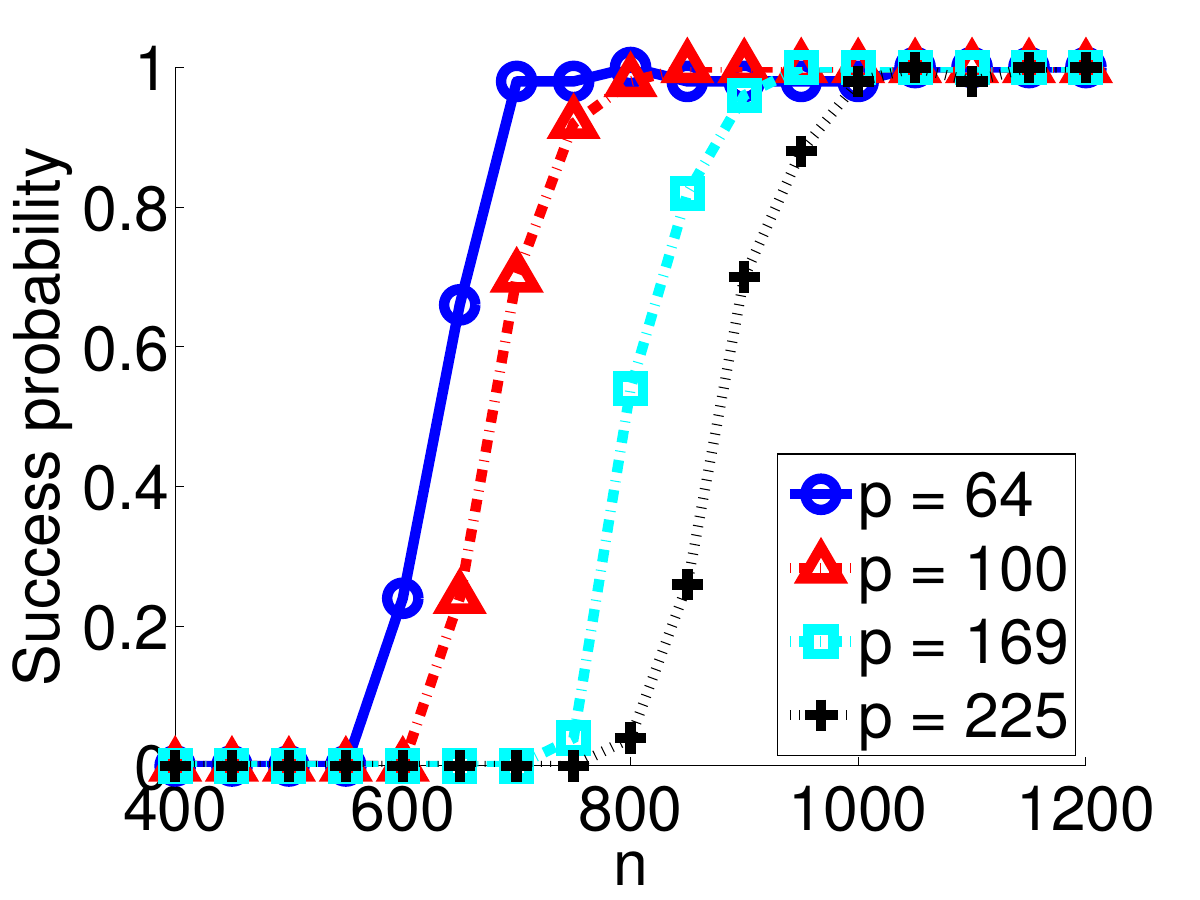}
        \includegraphics[width=0.4\textwidth]{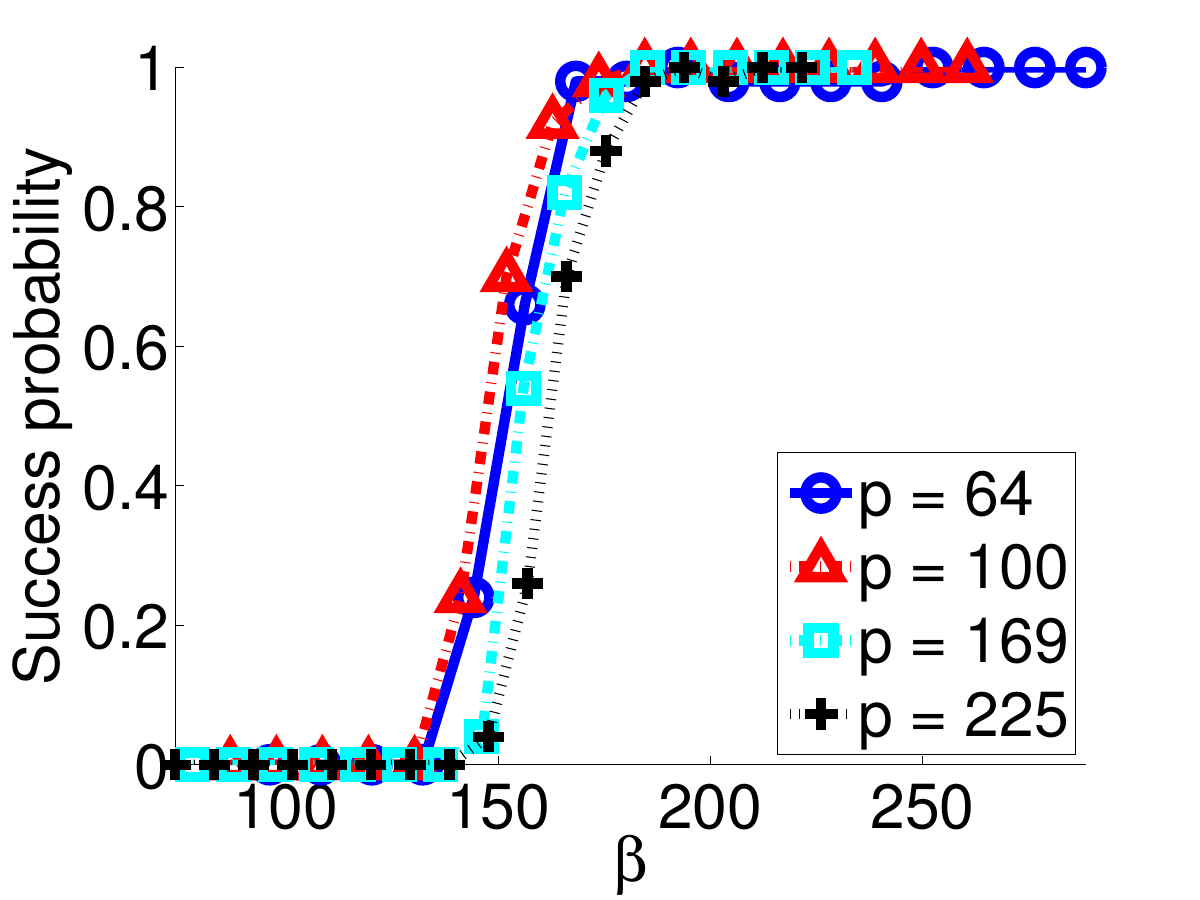}
    }
    \caption{Probabilities of successful support recovery for the (a)
    exponential MRF, grid structure with parameters $\theta^*_r =
    0.1$ and $\theta^*_{rt} =1$,  and the (b) Poisson MRF, grid structure
    with parameters $\theta^*_r = 2$ and $\theta^*_{rt} =-0.1$. The
    empirical
    probability of successful edge recovery over 50 replicates is
    shown versus the sample size $n$
    (left), and verses the re-scaled sample size $\beta = n/( \log
    p)$ (right).  The empirical curves align for the latter, thus
    verifying the logarithmic dependence of $n$ on $p$ as obtained in our sparsistency analysis. }
\end{figure}

\subsection{Simulation Studies}
We provide a small simulation study that corroborates our sparsistency
results; specifically Corollary~\ref{Cor:Ising} for
the exponential graphical model, where node-conditional distributions
follow an exponential distribution, and Corollary~\ref{Cor:Poisson}
for the
Poisson graphical model, where node-conditional distributions follow a
Poisson distribution. We instantiated the corresponding exponential
and Poisson graphical model distributions in \eqref{EqnExp} and
\eqref{EqnPoisson} for 4 nearest neighbor lattice graphs ($d=4$),  with varying number of nodes, $p \in
\{64,100,169,225\}$, and with identical edge weights for all edges:
 for exponential MRF, $\theta^*_r = 0.1$ and $\theta^*_{rt} =1$,
and, for Poisson MRF, $\theta^*_r = 2$ and $\theta^*_{rt}
=-0.1$. We generated i.i.d. samples from these distributions using
Gibbs sampling, and solved our sparsity-constrained $M$-estimation
problem by setting $\lambda_n = c \sqrt{\frac{\log p}{n}}$, following
our
corollaries; $c = 3$ for exponential MRF, and 15 for Poisson MRF. We
repeated each simulation 50 times and measured the empirical
probability over the 50 trials that our penalized graph estimate in
\eqref{EqnGraphEstimate} successfully recovered all edges, that is,
$P(\hat{E} = E^*)$. The left panels of
Figure~\ref{Fig:4nn_exponential} and
Figure~\ref{Fig:4nn_poisson} show the empirical probability of
successful edge recovery. In the right panel, we plot the empirical
probability against a re-scaled sample size $\beta = n/(\log
p)$. According to our corollaries, the sample size $n$ required for
successful graph structure recovery scales logarithmically with the
number of nodes $p$. Thus, we would expect the
empirical curves for different problem sizes to more closely align
with this re-scaled sample size on the horizontal axis, a result
clearly seen in the right panels
of Figure 1. This small numerical study thus corroborates our
theoretical sparsistency results.

\begin{figure}[h!tb]
    \centering
    \subfigure[Exponential Grid Structure]{\label{Fig:ROC_exponential}
        \centering
        \includegraphics[width=0.4\textwidth]{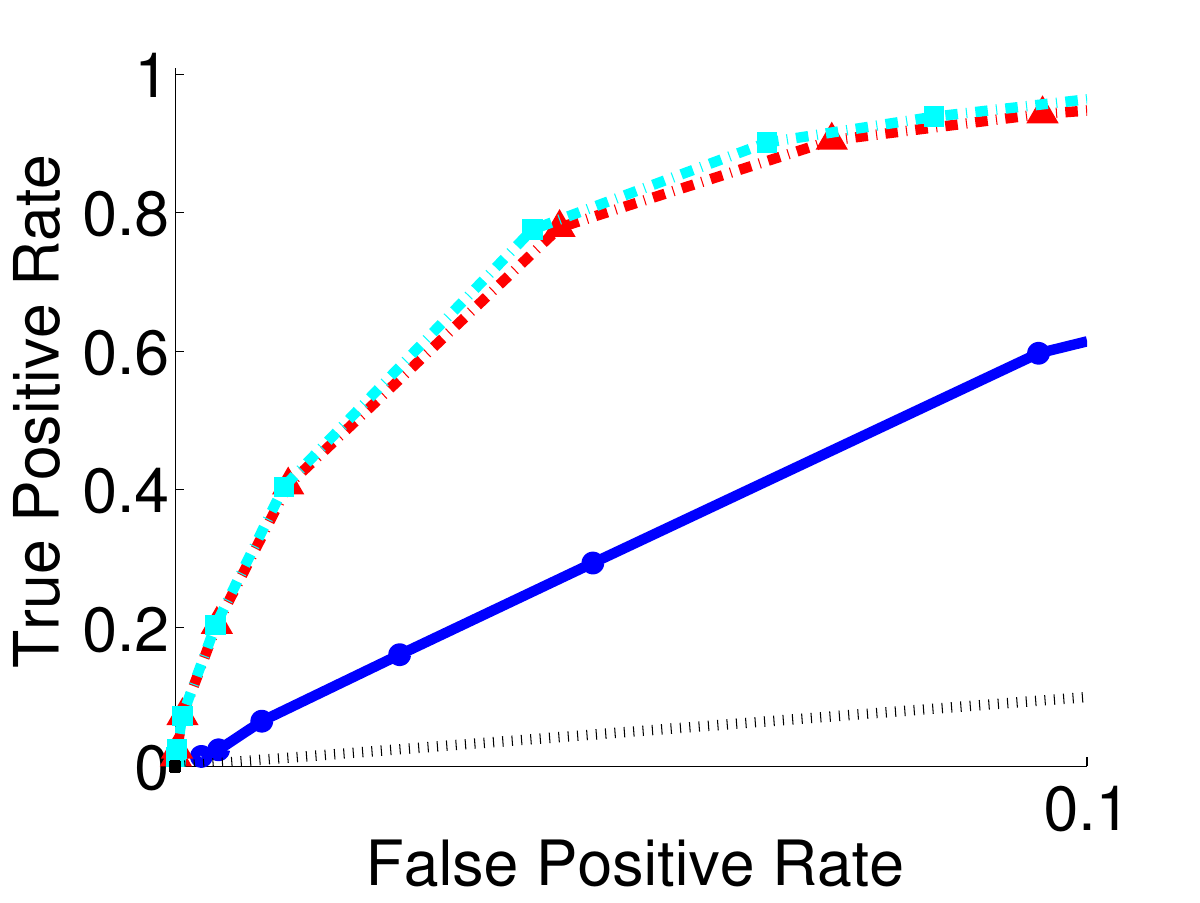}
    }
    \subfigure[Poisson Grid Structure]{\label{Fig:ROC_poisson}
        \centering
        \includegraphics[width=0.4\textwidth]{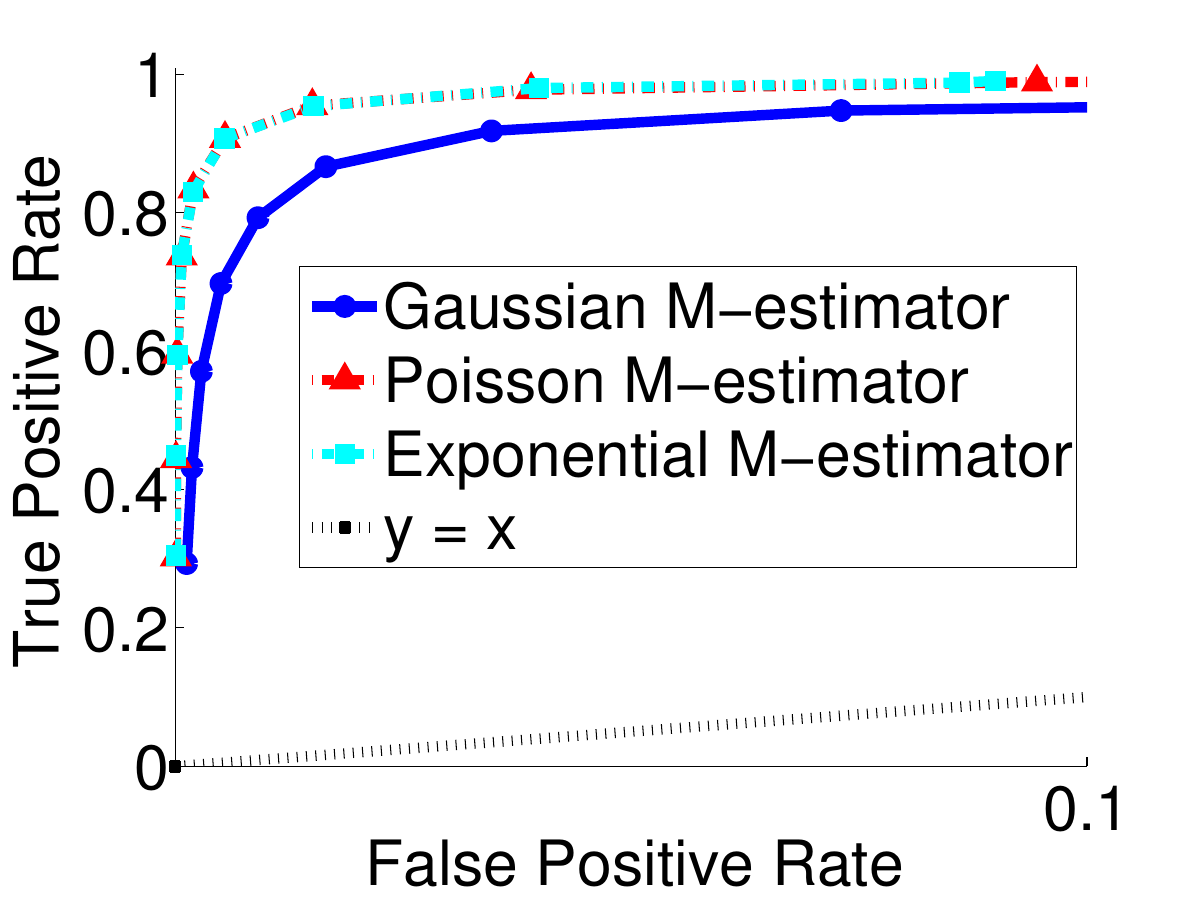}
    }
    \caption{Receiver-operator curves (ROC) computed by varying the
    regularization parameter, $\lambda_{n}$.  High-dimensional data is
    generated
    according to (a) the Exponential MRF with $(n,p) = (150,225)$ and to
    (b) the Poisson MRF with $(n,p) = (100,225)$.  Results are
    compared for three $M$-estimators: that of the Poisson,
    exponential, and Gaussian distributions. }
    \label{Fig:ROC_neg}
\end{figure}

We also evaluate the comparative performance of our $M$-estimators for
recovering the true edge structure from the different types of
data. Specifically, we consider the three typical examples in our
unified neighborhood selection approach: the Poisson $M$-estimator, the
Exponential $M$-estimator, and the well-known Gaussian
$M$-estimator by \citep{Meinshausen06}. In order to extensively
compare their performances, we compute the receiver-operator-curves
for the overall graph recovery by varying the regularization
parameter, $\lambda_{n}$. In Figure \ref{Fig:ROC_neg}, the same graph
structures for the exponential MRF ($\theta^*_r = 0.1$ and $\theta^*_{rt}
=1$) and the Poisson MRF ($\theta^*_r = 2$ and $\theta^*_{rt} =-0.1$) with
4 nearest neighbors, are used as in the previous simulation. Moreover,
we focus on the high-dimensional regime where $n < p $. As shown in the
figure, exponential and Poisson $M$-estimators outperform and have
significant advantage over Gaussian neighborhood selection approach if
the data is generated according to exponential or Poisson MRFs. One
interesting phenomenon we observe is that exponential and Poisson
$M$-estimators perform similarly regardless of the underlying graphical
model distribution.  This likely occurs as our estimator maximizes the
conditional likelihoods by fitting penalized GLMs.  Note that GLMs
assume that the conditional mean of the regression model follows an
exponential family distribution.  As both the Poisson distribution and
the exponential distribution have the same mean, the rate parameter,
$\lambda$, we would expect GLM-based methods that fit conditional
means to perform similarly.

\begin{figure}[h!tb]
    \centering
%    \subfigure[Gaussian copula modeling positively and negatively correlated exponential distributions]{\label{Fig:ROC_exponential_gcop}
%        \centering
%        \includegraphics[width=0.4\textwidth]{plots/plot_exponential_gcop_225_200.pdf}
%        \includegraphics[width=0.4\textwidth]{plots/plot_exponential_gcop_225_200negative.pdf}
%    }
    \subfigure{%[Gaussian copula modeling positively and negatively correlated Poisson  distributions]{\label{Fig:ROC_poisson_gcop}
        \centering
        \includegraphics[width=0.4\textwidth]{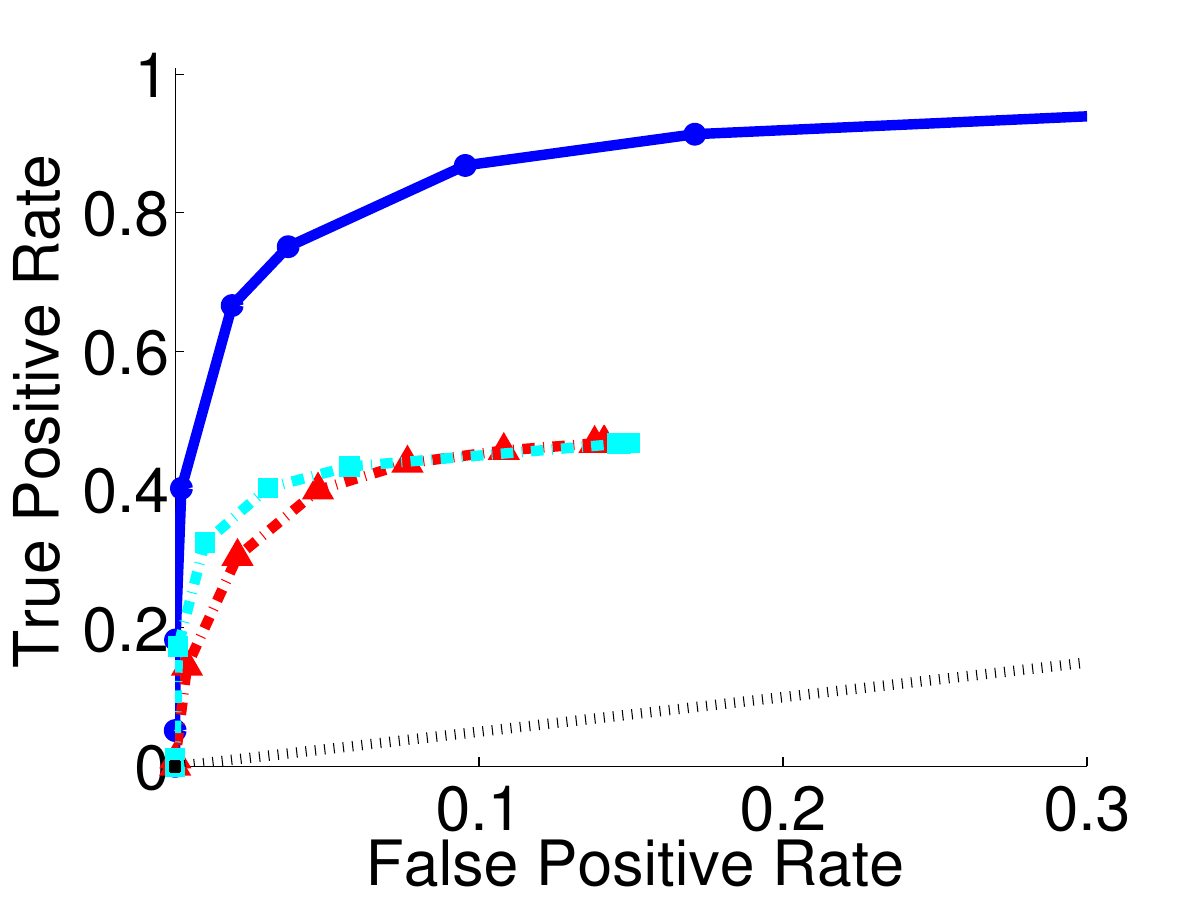}
        \includegraphics[width=0.4\textwidth]{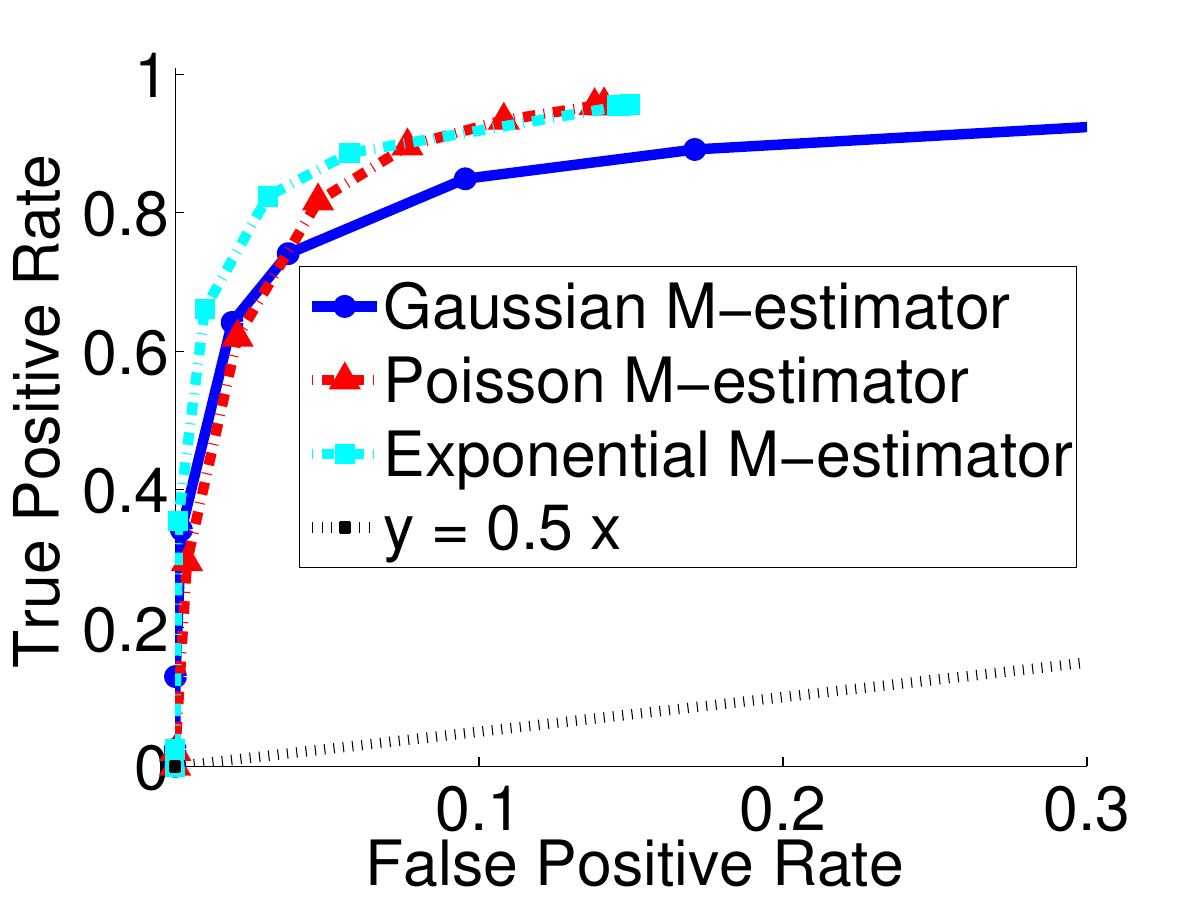}
    }
    \caption{Receiver-operator curves (ROC) computed by varying the
    regularization parameter, $\lambda_{n}$, for data, $(n,p) =
    (200,225)$, generated via Poisson
    copula transform according to a network with both positive and
    negative conditional dependencies. Left plot denotes results on overall edge recovery, while
    right plot denotes recovery of the edges with negative weights
    corresponding to negative conditional dependencies.  }
    \label{Fig:ROC_neg_cop}
\end{figure}

As discussed at end of Section \ref{gmodel}, the exponential and
Poisson graphical models are able to capture only negative conditional
dependencies between random variables, and our corresponding
$M$-estimators are computed under this constraint. In our last
simulation, we evaluate the impact of
this restriction when the true graph contains both
positive and negative edge weights.  As there does not exist a proper
MRF related to the Poisson and exponential distributions with both
positive and negative dependencies, we resort to generating data from
via a copula transform.  In particular, we first generate multivariate
Gaussian samples from $N(0, \Sigma$) where $\Theta = \Sigma^{-1}$ is the
precision matrix corresponding to the 4 nearest neighbor grid
structure previously considered.  Specifically, $\Theta$ has all ones
on the diagonal and $\theta_{rt}^{*} = \pm 0.2$ with equal
probabilities.  We then use a standard copula transform to make the
marginals of the generated data approximately Poisson.
Figure \ref{Fig:ROC_neg_cop} again present receiver operator curves (ROC) for the three different
classes of $M$-estimators on the copula transformed data, transformed to the
Poisson distribution. In the left of Figure~\ref{Fig:ROC_neg_cop}, we consider \emph{signed} support
recovery where we define the true positive rate as $\displaystyle
\frac{\text{$\#$ of edges s.t. $\text{sign}(\theta_{rt}^*) =
\text{sign}(\widehat{\theta}_{rt})$ }}{ \text{$\#$ of edges}  }$. In the
right, on the other hand, we \emph{ignore} the positive edges
so that true positive rate is now $\displaystyle \frac{\text{$\#$ of
edges s.t. $\text{sign}(\theta_{rt}^*) = \text{sign}(\widehat{\theta}_{rt}) = -1 $
}}{ \text{$\#$ of negative edges}  }$. Note that the false positive rate
is also defined similarly.  As expected, the results indicate that our
Poisson and
exponential $M$-estimators fail to recover the edges with positive
conditional dependencies recovered by the Gaussian $M$-estimator.
However, when attention is restricted to negative conditional
dependencies, our method outperforms the Gaussian $M$-estimator.
Notice also that for the exponential and Poisson
$M$-estimators, the highest false positive rate achieved is around
$0.15$.  This likely occurs due to the constraints enforced by our
$M$-estimators that force the weights of potential positive
conditional dependent edges to be zero.  Thus, while the restrictions
on the edge weights may be severe, for the purpose of estimating
negative conditional dependencies with limited false positives, the
Poisson and exponential $M$-estimators have an advantage.

\subsection{Real Data Examples}

%% \begin{figure}[!b]
%% \begin{center}
%% \includegraphics[width=3in,clip=true,trim=2.2in
%%   6.2in 2in .05in]{plots/miRNA-meta.pdf}\includegraphics[width=3in]{plots/fig_expgm_sachs.pdf}
%% \end{center}
%% \caption{Meta-miRNA inhibitory network for breast
%%   cancer learned via Poisson graphical models (left) and a inhibitory
%%   cell signaling network learned from flow cytometry data via
%%   exponential graphical models (right).   For the meta-miRNA inhibitory network,
%%   miRNA-sequencing data from TCGA was processed into tightly
%%   correlated clusters, meta-miRNAs, with the driver miRNAs identified
%%   for each cluster taken as the set of nodes for our network.  The
%%   Poisson network reveals major inhibitory
%%   relationships between three hub miRNAs, two of which have been
%%   previously identified as tumor suppressors in breast cancer.  For
%%   the cell signaling network, an exponential graphical model was fit
%%   to un-transformed flow cytometry data measuring 11 proteins.  The
%%   exponential network identifies PKA (protein kinase A) as a major
%%   inhibitor, consistent with previous results.}
%% \label{fig_example}
%% \end{figure}

To demonstrate the versatility of our family of graphical models, we
also provide two real data examples: a meta-miRNA inhibitory
network estimated by the Poisson graphical model, Figure
\ref{fig_miRNA}, and a cell signaling network
estimated by the exponential graphical model, Figure~\ref{fig_expgm}.

When applying our family of graphical models, there is
always a question of whether our model is an appropriate fit for the observed
data.  Typically, one can assess model fit using goodness-of-fit
tests.  For the Gaussian graphical model, this reduces to testing
whether the data follows a multivariate Gaussian distribution.  For
general exponential family graphical models, testing for
goodness-of-fit is more challenging.  Some have proposed likelihood
ratio tests specifically for lattice systems with a fixed and known
dependence structure \citep{besag_1974}.  When the network structure
is unknown,
however, there are no such existing
tests.  While we leave the development of an exact test to future
work, we provide a heuristic that can help us understand whether
our model is appropriate for a given dataset.

Recall that our model
assumes that conditional on its node-neighbors, each variable is
distributed according to an exponential family.  Thus, if the
neighborhood is known, our conditional models are simply GLMs, for which the goodness-of-fit can be assessed compared
to a null model by a likelihood ratio test
\citep{McCullagh}.  When
neighborhoods must be estimated, and specifically when estimated via
an $\ell_{1}$-norm penalty, the resulting ratio of likelihoods no
longer follow a chi-squared distribution
\citep{buhlmann2011statistics}.  Recently, for the
$\ell_{1}$ linear regression case, \citet{lockhart2013significance}
have shown that the
difference in the residual sums of squares follows an exponential
distribution.  Similar results have not yet been extended to the
penalized GLM case.  In the absence of such tests, we propose a simple
heuristic: for each node, first estimate the node-neighborhood via our
proposed M-estimator.  Next, assuming the neighborhood is fixed, fit a
GLM and compare the fit of this model to that of a null model (only an
intercept term) via the likelihood ratio test.  One can then
heuristically assess the overall goodness-of-fit by examining the fit
of a GLM to all the nodes.  This procedure is clearly not an exact
test, and following from \citet{lockhart2013significance}, it is
likely conservative.  In the
absence of an exact test, which we leave for future work, this
heuristic provides some assurances about the appropriateness of our
model for real data.

\subsubsection{Poisson Graphical Model: Meta-miRNA Inhibitory Network}

\begin{figure}[h!b]
\begin{center}
\centering
\includegraphics[width=4in,clip=true,trim=3in 1.5in 5in 0.5in]{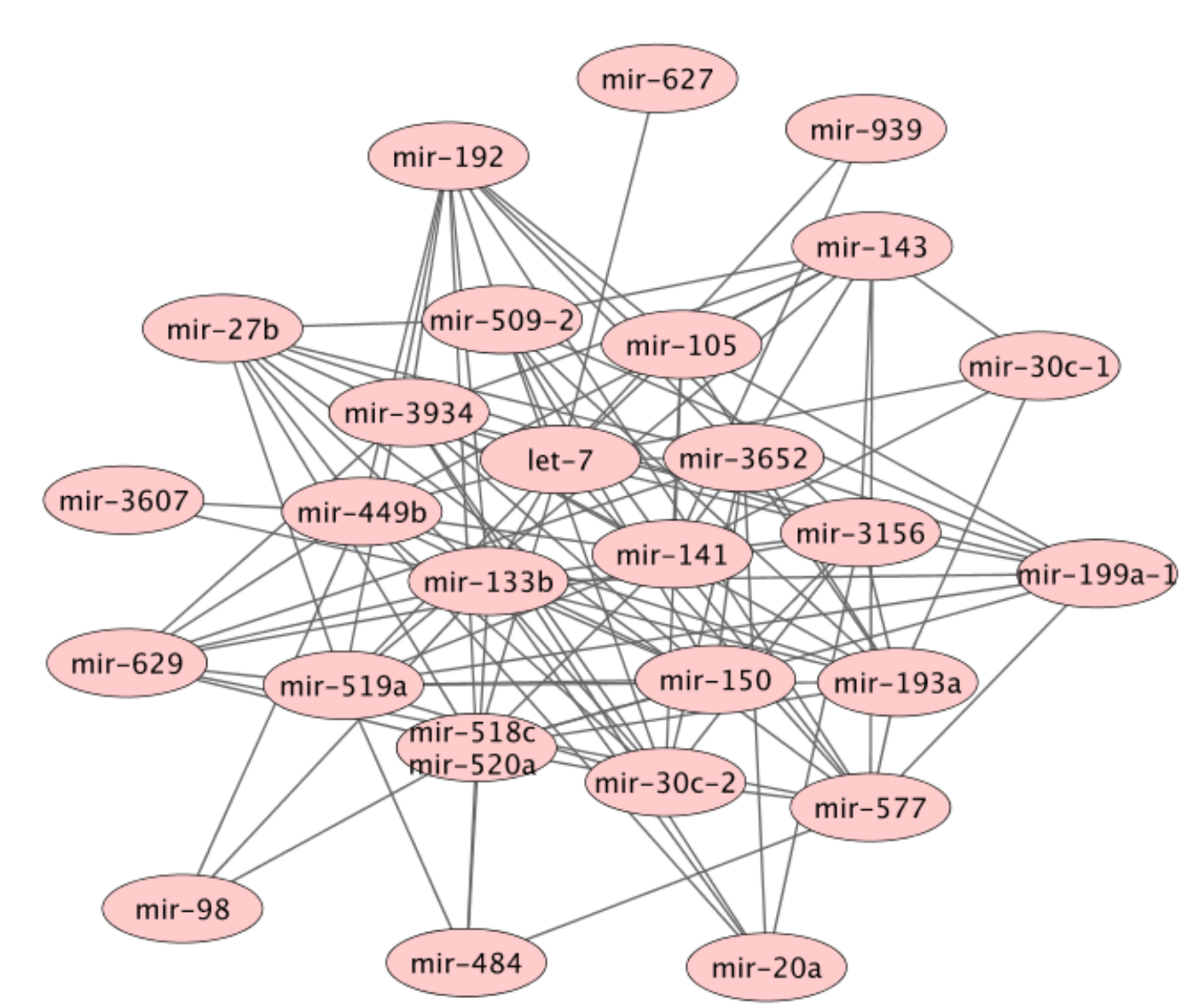}
\end{center}
\caption{Meta-miRNA inhibitory network for breast
  cancer estimated via Poisson graphical models from miRNA-sequencing
  data.   Level III data from TCGA was processed into tightly
  correlated clusters, meta-miRNAs, with the driver miRNAs identified
  for each cluster taken as the set of nodes for our network.  The
  Poisson network reveals major inhibitory
  relationships between three hub miRNAs, two of which have been
  previously identified as tumor suppressors in breast cancer. }
\label{fig_miRNA}
\end{figure}

Gaussian graphical models have often been used to study
high-throughput genomic networks estimated from microarray data
\citep{Peer:2001ci,Friedman:2004uu,wei2007markov}. Many
high-throughput technologies, however, do not produce even approximately
Gaussian data, so that our class of graphical models could be
particularly important for estimating genomic networks from
such data. We demonstrate the applicability of our class of models by
estimating a meta-miRNA inhibitory network for breast cancer estimated
by a Poisson graphical model. Level III breast cancer miRNA expression
\citep{tcga_breast_2012} as measured by next generation sequencing was
downloaded from the TCGA  portal
(http://tcga-data.nci.nih.gov/tcga/). MicroRNAs (miRNA) are short RNA
fragments that are thought to be post-transcriptional regulators,
predominantly inhibiting translation.  Measuring miRNA expression by
high-throughput sequencing results in \emph{count data} that is
zero-inflated, highly skewed, and whose total count volume depends on
experimental conditions \citep{li_rna_2011}. Data was processed to be
approximately Poisson by following the steps described in
\citep{allen2013local}. In brief, the data was quantile corrected to
adjust for sequencing depth \citep{bullard_rna_2010}; the miRNAs with
little variation across the samples, the bottom 50\%, were filtered
out; and the data was adjusted for possible over-dispersion using a
power transform and a goodness of fit test
\citep{li_rna_2011}. We also tested for batch effects
in the resulting data matrix consisting of 544 subjects and 262
miRNAs: we fit a Poisson ANOVA model \citep{leek2010tackling}, and
only found 4\% of miRNAs to be associated with batch labels; and thus
no significant batch association was detected.  As several miRNAs
likely target the same gene or genes in the same pathway, we expect
there to be strong positive dependencies among variables that cannot
be captured directly by our Poisson graphical model which only permits
negative conditional relationships.  Thus, we will use our model to
study inhibitory relationships between what we term meta-miRNAs, or
groups of miRNAs that are tightly positively correlated.
To accomplish this, we further processed our data to form
clusters of positively correlated miRNAs using hierarchical clustering
with average linkage and one minus the correlation as the distance
metric.  This resulted in 40 clusters of tightly positively correlated
miRNAs.  The nodes of our meta-miRNA network were then taken as a the
medoid, or median centroid defined as the miRNA closest in Euclidean
distance to the cluster centroid, in each group.

A Poisson graphical model was fit to the meta-miRNA data by performing
neighborhood selection with the sparsity of the graph determined by
stability selection \citep{liu_STARS_2010}.  The heuristic previously
discussed was used to assess goodness-of-fit for our model.  Out of
the 40 node-neighborhoods tested via a likelihood ratio test, 36
exhibited $p$-values less than 0.05, and 34 were less than 0.05/40,
the Bonferroni-adjusted significance level.  These results show that the
Poisson GLM is a significantly better fit for the majority of
node-neighborhoods than the null model, indicating that our Poisson
graphical model is appropriate for this data.  The results of our
estimated Poisson graphical model,
Figure~\ref{fig_miRNA} (left), are consistent with the cancer genomics
literature.  First, the meta-miRNA inhibitory network has three major
hubs.  Two of these, miR-519 and miR-520, are known to be breast
cancer tumor suppressors, suppressing growth by reducing HuR
levels~\citep{abdelmohsen2010mir} and by targeting NF-KB and TGF-beta
pathways \citep{keklikoglou2011microrna} respectively. The third major
hub, miR-3156, is a miRNA of unknown function; from its major role in
our network, we hypothesize that miR-3156 is also associated with
tumor suppression.  Also interestingly, let-7, a well-known miRNA
involved in tumor metastasis \citep{yu2007let}, plays a central role
in our network, sharing edges with the five largest hubs.  This
suggests that our Poisson graphical model has recovered relevant
negative relationships between miRNAs with the five major hubs acting
as suppressors, and the central let-7 miRNA and those connected to
each of the major hubs acting as enhancers of tumor
progression in breast cancer.

\subsubsection{Exponential Graphical Model: Inhibitory Cell-Signaling Network}

\begin{figure}[h!b]
\begin{center}
\includegraphics[width=6in]{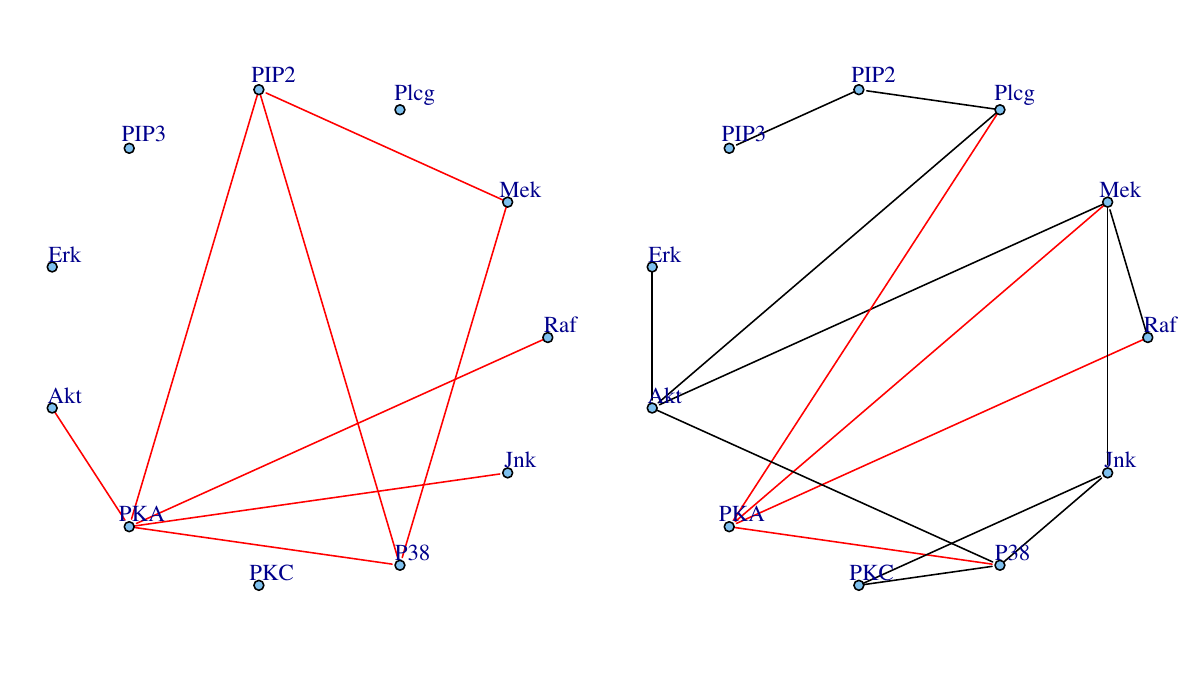}
\end{center}
\caption{ Cell signaling network estimated from flow
cytometry data via exponential graphical models (left) and Gaussian
  graphical models (right).   The exponential graphical model was fit
  to un-transformed flow cytometry data measuring 11 proteins, and the
  Gaussian graphical model to log-transformed data.  Estimated
  negative conditional dependencies are given in red.
  Both networks identify PKA (protein kinase A) as a major
  inhibitor, consistent with previous results.}
\label{fig_expgm}
\end{figure}

We demonstrate our exponential graphical model, derived from the
univariate exponential distribution, using a protein signaling example
\citep{sachs2005causal}.  Multi-florescent flow cytometry
was used to measure the presence of eleven proteins ($p=11$) in
$n=7462$ cells.  This data set was first analyzed using Bayesian
Networks in \citet{sachs2005causal} and then using the graphical lasso
algorithm in \citet{glasso}.  Measurements from flow-cytometry data
typically follow a left skewed distribution. Thus to model such data,
these measurements are typically normalized to be approximately
Gaussian using a log transform after shifting the data to be
non-negative \citep{herzenberg2006interpreting}. Here, we demonstrate
the applicability of our exponential graphical models to recover
networks directly from continuous skewed data, so that we learn the
network directly from the flow-cytometry data without any log or such
transforms.  Our pre-processing is limited to shifting the data for
each protein so that it consists of non-negative values.  For
comparison purposes, we also fit a Gaussian graphical model to the
log-transformed data.

We then learned an exponential and Gaussian graphical model from this
flow cytometry
data using stability selection \citep{liu_STARS_2010} to select the
sparsity of the graphs.  The goodness-of-fit heuristic previously
described was used to assess the appropriateness of our model.  Out of
the eight connected node-neighborhoods, the likelihood ratio test was
statistically significant for seven neighborhoods, indicating that our
exponential GLM is a better fit than the null model.
The estimated protein-signaling network is
shown on the right in Figure~\ref{fig_expgm} with that of the Gaussian
graphical model fit to the log-transformed data on the left. Estimated
negative conditional dependencies are shown in red.  Recall that the
exponential graphical model restricts the edge weights to be
non-negative; because of the negative inverse link, this implies that
only negative conditional associations can be estimated.  Notice that
our exponential graphical model finds that PKA, protein kinase A, is a
major protein
inhibitor in cell signaling networks.  This is consistent with the
inhibitory relationship of PKA as estimated by the Gaussian graphical
model, right Figure~\ref{fig_expgm}, as well as its hub status in the
Bayesian network of \citep{sachs2005causal}.  Interestingly, our
exponential graphical
model also finds a clique between PIP2, Mek, and P38, which was not
found by Gaussian graphical models.

%%%%%%%%%%%%%%%%%%%%%%%%%%%%%%%%%%%%%%%%%%%%%%%%%%%%%%%
\section{Discussion}

We study what we call the class of exponential family graphical models that arise when we assume that node-wise conditional distributions follow exponential family distributions.  Our work broadens the class of off-the-shelf graphical models from classical instances such as Ising and Gaussian graphical models.  In particular, our class of graphical models provide closed form multivariate densities as extensions of several univariate exponential family distributions (e.g. Poisson, exponential, negative binomial) where few currently exist; and thus may be of further interest to the statistical community. Further, we provide simple M-estimators for estimating any of these graphical models from data, by fitting node-wise penalized conditional exponential family distributions, and show that these estimators enjoy strong statistical guarantees. The statistical analyses of our M-estimators required subtle techniques that may be of general interest in the analysis of sparse M-estimation.

There are many avenues of future work related to our proposed models. We assume that all conditional distributions are members of an exponential family.  To determine whether this assumption is appropriate in practice for real data, a goodness-of-fit procedure is needed.  While we have proposed a heuristic to this effect, more work is needed to determine a rigorous likelihood ratio test for testing model fit. For several instances of our proposed class of models, specifically those with variables with infinite domains, severe restrictions on the parameter space are sometimes needed.  For instance, the Poisson and exponential graphical models studied in Section 4, could only model negative conditional dependencies, which may not always be desirable in practice.  A key question for future work is whether these
restrictions can be relaxed for particular exponential family distributions. Finally, while we have focused on single parameter exponential families, it would be interesting to investigate the consequences of using multi-parameter exponential family distributions. Overall, our work has opened avenues for learning Markov networks from a broad class of univariate distributions, the properties and applications of which leave much room for future research.

% Acknowledgements should go at the end, before appendices and references

\acks{We would like to acknowledge support for this project from ARO
  W911NF-12-1-0390, NSF IIS-1149803, IIS-1320894, IIS-1447574, and
  DMS-1264033 (PR and EY); NSF DMS-1264058 and DMS-1209017 (GA); and
  the Houston Endowment and NSF DMS-1263932 (ZL).}

% Manual newpage inserted to improve layout of sample file - not
% needed in general before appendices/bibliography.

%\newpage

%%%%%%%%%%%%%%%%%%%%%%%%%%%%%%%%%%%%%%%%%%%%%%%%%%%%%%%
\appendix
\def\bxs{\overline{X}_{r:0}}
\def\bxi{\overline{X}_{1:0}}

\section{Proof of Theorem \ref{PropGeneralCondForm}}

The proof follows the development in \cite{besag_1974}, where they consider the case with $k=2$. We define $Q(X)$ as $Q(X) := \log(P(X)/P({\bf 0}))$, for any $X = (X_1,\hdots,X_p) \in \mathcal{X}^{p}$ where $P({\bf 0}))$ denotes the probability that all random variables take 0. Given any $X$, also denote $\bxs := (X_1,\hdots,X_{r-1},0,X_{r+1},\hdots,X_p)$. Now, consider the following the most general form for $Q(X)$:
\begin{align}\label{EqnQFactorForm}
	Q(X) = \sum_{1\leq r \leq p} X_{r} G_{r}(X_{r}) + \hdots + \sum_{1 \leq r_1 < r_2 < \hdots < r_k \leq p} X_{r_1}\hdots X_{r_k} G_{r_1 \hdots r_k}(X_{r_1},\hdots,X_{r_k}),
\end{align}
since the joint distribution has factors of at most size $k$. By the definition of $Q$ and some algebra (See Section 2 of \cite{besag_1974} for details), it can then be seen that
\begin{align}\label{EqnTmpA}
	& \exp(Q(X) - Q(\bxs)) \nonumber\\
	= \ & P(X_r|X_1,\hdots,X_{r-1},X_{r+1},\hdots,X_p)/P(0|X_1,\hdots,X_{r-1},X_{r+1},\hdots,X_p).
\end{align}
Now, consider the simplifications of both sides of \eqref{EqnTmpA}. For notational simplicity, we fix $r=1$ for a while. Given the form of $Q(X)$ in \eqref{EqnQFactorForm}, we have
\begin{align}\label{EqnTmpB}
	Q(X) - Q(\bxi) = \ & X_{1}\bigg(G_1(X_1) + \sum_{2\leq t \leq p} X_t G_{1t}(X_1, X_t) + \hdots  \nonumber\\
	&+ \sum_{2 \leq t_2 < t_3 < \hdots < t_k \leq p} X_{t_2}\hdots X_{t_k} G_{1 t_2 \hdots t_k}(X_1,X_{t_2}\hdots,X_{t_k})\bigg).
\end{align}
By given the exponential family form of the node-conditional distribution specified in the statement, right-hand side of \eqref{EqnTmpA} become
\begin{align}\label{EqnTmpC}
	\hspace{-.2cm}\log \frac{P(X_1|X_2,\hdots,X_p)}{P(0|X_2,\hdots,X_p)} = E(X_{V \backslash 1})(B(X_1) - B(0)) + (C(X_1) - C(0)).
\end{align}
Setting $X_{t} = 0$ for all $t \neq 1$ in \eqref{EqnTmpB}
and \eqref{EqnTmpC}, we obtain,
\begin{align*}
	X_{1} G_1(X_1) = E({\bf 0})(B(X_1) - B(0)) + (C(X_1) - C(0)).
\end{align*}
Setting $X_{t_2} = 0$ for all $t_2 \not\in \{1,t\}$,
\begin{align*}
	X_{1} G_1(X_1) + X_{1} X_{t} G_{1t}(X_1,X_t) = E(0,\hdots,X_t,\hdots,0)(B(X_1) - B(0)) + (C(X_1) - C(0)).
\end{align*}
Recovering the index 1 back to $r$ yields
\begin{align*}
	&X_{r} G_r(X_r) = E({\bf 0})(B(X_r) - B(0)) + (C(X_r) - C(0)), \\
	&X_{r} G_r(X_r) + X_{r} X_{t} G_{rt}(X_r,X_t) = E(0,\hdots,X_t,\hdots,0)(B(X_r) - B(0)) + (C(X_r) - C(0)).
\end{align*}
Similarly,
\begin{align}
	X_{t} G_t(X_t) + X_{r} X_{t} G_{rt}(X_r,X_t) = E(0,\hdots,X_r,\hdots,0)(B(X_t) - B(0)) + (C(X_t) - C(0)).
\end{align}
From the above three equations, we obtain:
\begin{align*}
	X_{r} X_{t} G_{rt}(X_r,X_t) = \theta_{rt} (B(X_r) - B(0)) (B(X_t) - B(0)).
\end{align*}
More generally, by considering non-zero triplets, and setting $X_{v} = 0$ for all $v \not\in \{r,t,u\}$, we obtain,
\begin{align}
	&X_{r} G_r(X_r) + X_{r} X_{t} G_{rt}(X_r,X_t) + X_{r} X_{u} G_{r u}(X_r,X_{u}) + X_{r} X_{t} X_{u} G_{rtu}(X_r,X_t,X_{u}) = \nonumber\\
    &E(0,\hdots,X_t,\hdots,X_{u},\hdots,0)(B(X_r) - B(0)) + (C(X_r) - C(0)),
\end{align}
so that by a similar reasoning we can obtain
\begin{align*}
	X_{r} X_{t} X_{u} G_{rtu}(X_r,X_t,X_{u}) = \theta_{rtu}(B(X_{r}) - B(0))(B(X_{t}) - B(0))(B(X_{u}) - B(0)).
\end{align*}
More generally, we can show that
\begin{align*}
	X_{t_1}\hdots X_{t_k} G_{t_1,\hdots,t_k}(X_{t_1},\hdots,X_{t_k}) = \theta_{t_1,\hdots,t_k}(B(X_{t_1}) - B(0))\hdots(B(X_{t_k}) - B(0)).
\end{align*}
Thus, the $k$-th order factors in the joint distribution as specified in \eqref{EqnQFactorForm} are tensor products of $(B(X_r) - B(0))$, thus proving the statement of the theorem.

\section{Proof of Proposition~\ref{PropSubGaussVariableGLM}}

By the definition of \eqref{ABar} with the following simple calculation, the moment generating function of $X_r^2$ becomes:
%Suppose we zero-pad the edge-weights in the true parameter $\theta^* \in \reals^{p+\binom{p}{2}}$ of \eqref{EqnJointGLM} to include a zero weight for sufficient statistic $X_{r}^2$; we will overload notation and denote this zero-padded parameter in $\reals^{p+\binom{p}{2} + 1}$ as $\theta^*$. Similarly, let $\bar{v} \in \reals^{p+\binom{p}{2} + 1}$ be the zero-padded parameter with its last coordinate equal to $t \in \reals$, so that $\|\bar{v}\|_2 = t$. A simple calculation shows that
%{\small
\begin{align*}
    \log \E[\exp(a X_r^2)] &= \log \int_{X} \exp\Big\{ a X_r^2 + \sum_{t \in V} \theta_t^* X_t + \sum_{(t,u) \in E }\theta_{tu}^* \, X_t \, X_{u} + \sum_{t \in V} C(X_t) - A(\theta^*) \Big\} \\
                           &= \bar{A}_r(a;\theta^*) - \bar{A}_r(0;\theta^*).
\end{align*}%}
Suppose that $a \leq 1$. Then, by a Taylor Series expansion, we have for some $\nu \in [0,1]$,
\begin{align*}
	\bar{A}_r(a;\theta^*) - \bar{A}_r(0;\theta^*) &= a \frac{\partial }{\partial\eta}\bar{A}_r(0;\theta^*) + \frac{1}{2} a^2 \frac{\partial^{2}}{\partial\eta^2}\bar{A}_r(\nu a;\theta^*)									\le \ConstVar a + \frac{1}{2}\ConstHess a^2,
\end{align*}
where the inequality uses Condition \ref{AN_GLM_MH}. Note that since the derivative of log-partition function is the mean of the corresponding sufficient statistics and $\bar{A}_r(0;\theta) = A(\theta)$, $\frac{\partial }{\partial\eta}\bar{A}_r(0;\theta^*) = \E[X_r^2] \leq \ConstVar$ by assumption.
Thus, by the standard Chernoff bounding technique, for all positive $a \leq 1$,
\begin{align*}
	P\bigg(\frac{1}{n}\sum_{i=1}^n\big(\xs\big)^2 \geq \delta \bigg) &\leq
	\exp(- n \delta a + n \ConstVar a + \frac{n}{2} \ConstHess a^2)
\end{align*}
With the choice of $a = \frac{\delta - \ConstVar}{\ConstHess} \leq 1$, we obtain
\begin{align*}
	P\bigg(\frac{1}{n}\sum_{i=1}^n\big(\xs\big)^2 \geq \delta \bigg) \le \exp\bigg( - n\frac{(\delta - \ConstVar)^2}{2 \ConstHess}\bigg) \le \exp\Big( - n\frac{\delta^2}{8 \ConstHess}\Big)
\end{align*}
provided that $\delta \leq 2\ConstVar/3$, as in the statement.

\section{Proof of Proposition~\ref{PropSubGaussSingleVariableGLM}}
Let $\bar{v} \in \reals^{p+\binom{p}{2}}$ be the zero-padded parameter with only one non-zero coordinate, which is $1$, for the sufficient statistics $X_r$ so that $\|\bar{v}\|_2 = 1$. A simple calculation shows that
\[
\log \E[\exp(X_r)] = A(\theta^* + \bar{v}) - A(\theta^*).
\]
By a Taylor Series expansion and Condition \ref{AN_GLM_MH}, we have for some $\nu \in [0,1]$,
\begin{align*}
A(\theta^* + \bar{v}) - A(\theta^*) &= \grad A(\theta^*) \cdot \bar{v} + \frac{1}{2}\bar{v}^{T} \grad^{2}A(\theta^* +  \nu \bar{v}) \bar{v}\\	
									&\overset{(i)}{\leq} \E[X_r] \|\bar{v}\|_2 + \frac{1}{2} \frac{\partial^{2}}{\partial\theta_r^2}A(\theta^* +  \nu \bar{v}) \|\bar{v}\|_2^2 \le \ConstMu + \frac{1}{2}\ConstHess ,
\end{align*}
where the inequality $(i)$ uses the fact that $\bar{v}$ has only nonzero element for the sufficient statistics $X_r$. Thus, again by the standard Chernoff bounding technique, for any positive constant $a$, $P\left(X_r \geq a \right) \leq \exp(- a + \ConstMu + \frac{1}{2} \ConstHess)$, and by setting $a = \delta \log \eta$ we have
\begin{align*}
	P\left(X_r \geq \delta \log \eta \right) &\leq \exp(- \delta \log \eta + \ConstMu + \frac{1}{2} \ConstHess) \leq c \eta^{-\delta} ,
\end{align*}
where $c = \exp(\ConstMu + \frac{1}{2}\ConstHess)$, as claimed.

\section{Proof of Theorem \ref{Thm:Main}}

In this section, we sketch the proof of Theorem \ref{Thm:Main} following the \emph{primal-dual witness} proof technique in \cite{Wainwright2006new,RWL10}. We first note that the optimality condition of the convex program~\eqref{Eq:Obj} can be written as
\begin{align}\label{Eq:subOptCond}
    \nabla \ell(\widehat{\theta};X^{1:n}) + \lambda_n \widehat{Z} = 0,
\end{align}
where $\widehat{Z}$ is a length $p$ vector: $\widehat{Z}_{\backslash r} \in \partial \|\widehat{\theta}_{\backslash r}\|_1$ is a length $(p-1)$ subgradient vector where $\widehat{Z}_{rt} = $ sign$(\widehat{\theta}_{rt})$ if $\widehat{\theta}_{rt} \neq 0$, and $|\widehat{Z}_{rt}| \leq 1$ otherwise; while $\widehat{Z}_r$, corresponding to $\theta_r$, is set to 0 since the nodewise term $\theta_r$ is not penalized in the $M$-estimation problem~\eqref{Eq:Obj}.

Note that in a high-dimensional regime with $p \gg n$, the convex program~\eqref{Eq:Obj} is not necessarily \emph{strictly} convex, so that it might have multiple optimal solutions. However, the following lemma, adapted from \cite{RWL10}, shows that nonetheless the solutions share their support set under certain conditions. We first recall the  notation $S = \{(r,t) : t \in N^*(r)\}$ to denote the true neighborhood of node $r$, and $S^c$ to denote its complement.
\begin{lemma}\label{LemPDTmp}
Suppose that there exists a primal optimal solution $\widehat{\theta}$ with associated subgradient $\widehat{Z}$ such that $\|\widehat{Z}_{S^c}\|_{\infty} < 1$. Then, any optimal solution $\tilde{\theta}$ will satisfy $\tilde{\theta}_{S^c} = 0$. Moreover, under the condition of $\|\widehat{Z}_{S^c}\|_{\infty} < 1$, if $Q_{SS}^*$ is invertible, then $\widehat{\theta}$ is the unique optimal solution of \eqref{Eq:Obj}.
\end{lemma}

\begin{proof}
This lemma can be proved by the same reasoning developed for the special cases \cite{Wainwright2006new,RWL10} in our framework; As in the previous works, for any node-conditional distribution in the form of exponential family, we are solving the convex objective with $\ell_1$ regularizer \eqref{Eq:Obj}. Therefore, the problem can be written as an equivalent constrained optimization problem, and by the complementary slackness, for any optimal solution $\tilde{\theta}$, we have $\inner{\widehat{Z}, \tilde{\theta}} = \| \tilde{\theta} \|_1$. This simply implies that for all index $j$ for which $ | \widehat{Z}_j| <1 $, $\tilde{\theta}_j =0$ (See \cite{RWL10} for details). Therefore, if there exists a primal optimal solution $\widehat{\theta}$ with associated subgradient $\widehat{Z}$ such that $\|\widehat{Z}_{S^c}\|_{\infty} < 1$, then, any optimal solution $\tilde{\theta}$ will satisfy $\tilde{\theta}_{S^c} = 0$ as claimed.

Finally, we consider the restricted optimization problem subject to the constraint $\theta_{S^C} = 0$. For this restricted optimization problem, if the Hessian, $Q_{SS}^*$, is positive definite as assumed in the lemma, then, this restricted problem is strictly convex, and its solution is unique. Moreover, since all primal optimal solutions of \eqref{Eq:Obj}, $\tilde{\theta}$, satisfy $\tilde{\theta}_{S^c} = 0$ as discussed, the solution of the restricted problem is the unique solution of \eqref{Eq:Obj}.
\end{proof}

We use this lemma to prove the theorem following the \emph{primal-dual witness} proof technique in \cite{Wainwright2006new,RWL10}. Specifically, we \emph{explicitly construct} a pair $(\widehat{\theta},\widehat{Z})$ as follows (denoting the true support set of the edge parameters by $S$):
\begin{enumerate}[leftmargin=0.05cm, itemindent=0.85cm,label=\textbf{(\alph*)}]
    \item Recall that $\theta(r) = (\theta_r, \theta_{\backslash r}) \in \reals\times\reals^{p-1}$. We first fix $\theta_{S^c} =0$ and solve the restricted optimization problem: $(\widehat{\theta}_r,\, \widehat{\theta}_S,0) = \argmin_{\theta_r \in \reals, \, (\theta_S,0)\in \reals^{p-1}} \{\ell(\theta;X^{1:n}) + \lambda_n \|\theta_S\|_1\}$, and $\widehat{Z}_S =$ sign$(\widehat{\theta}_S)$.
    \item We set $\widehat{\theta}_{S^c} = 0$.
    \item We set $\widehat{Z}_{S^c}$ to satisfy the condition \eqref{Eq:subOptCond} with $\widehat{\theta}$ and $\widehat{Z}_S$.
\end{enumerate}

By construction, the support of $\widehat{\theta}$ is included in the true support $S$ of $\theta^*$, so that we would finish the proof of the theorem provided (a) $\widehat{\theta}$ satisfies the stationary condition of ~\eqref{Eq:Obj}, as well as the condition $\|\widehat{Z}_{S^c}\|_{\infty} < 1$ in Lemma~\ref{LemPDTmp} with high probability, so that by Lemma~\ref{LemPDTmp}, the primal solution $\widehat{\theta}$ is guaranteed to be unique; and (b) the support of $\widehat{\theta}$ is not strictly within the true support $S$. We term these conditions \emph{strict dual feasibility}, and \emph{sign consistency} respectively.

We will now rewrite the subgradient optimality condition~\eqref{Eq:subOptCond} as
$$\nabla^2 \ell(\theta^*;X^{1:n})(\widehat{\theta}-\theta^*) = - \lambda_n \widehat{Z} + W^n + R^n,$$
where $W^n := - \nabla\ell(\theta^*;X^{1:n})$ is the sample score function (that we will show is small with high probability), and $R^n$ is the remainder term after coordinate-wise applications of the mean value theorem; $R_j^n = [\nabla^2 \ell(\theta^*;X^{1:n})- \nabla^2 \ell(\bar{\theta}^{(j)};X^{1:n})]_j^T(\widehat{\theta}-\theta^*)$, for some $\bar{\theta}^{(j)}$ on the line between $\widehat{\theta}$ and $\theta^*$, and with $[\cdot]_j^T$ being the $j$-th row of a matrix.

Recalling the notation for the Fisher information matrix $Q^* := \nabla^2 \ell(\theta^*;X^{1:n})$, we then have
\[ Q^* (\widehat{\theta}-\theta^*) = - \lambda_n \widehat{Z} + W^n + R^n.\]

From now on, we provide lemmas that respectively control various terms in the above expression: the score term $W^n$, the deviation $\widehat{\theta} -\theta^*$, and the remainder term $R^n$.
The first lemma controls the score term $W^n$:
\begin{lemma}\label{Lem:1}
    Suppose that we set $\lambda_n$ to satisfy $\frac{8(2-\alpha)}{\alpha}\sqrt{\DConstOne\ConstGLM}\sqrt{\frac{\log p}{n}} \leq \lambda_n \leq  \frac{4(2-\alpha)}{\alpha} \allowbreak \DConstOne \DConstTwo \ConstGLM$ for some constant $\ConstGLM \leq \min\{2\ConstVar/3,2\ConstHess+\ConstVar\}$. Suppose also that $n \geq \frac{8\ConstHess^2}{\ConstGLM^2}\log p$. Then, given a incoherence parameter $\alpha \in (0,1]$,
    \begin{align*}
        P\left(\frac{2-\alpha}{\lambda_n}\|W^n\|_{\infty} \leq \frac{\alpha}{4}\right) \geq 1- c_1 p'^{-2} - \exp(-c_2 n) - \exp(-c_3 n)
    \end{align*}
    where $p' := \max\{n,p\}$.
\end{lemma}
\noindent The next lemma controls the deviation $\widehat{\theta} -\theta^*$:
\begin{lemma}\label{Lem:2}
    Suppose that $\lambda_n d \leq \frac{\Cmin^2}{40 \Dmax \DConstThree \log p'}$ and $\|W^n\|_{\infty} \leq \frac{\lambda_n}{4}$. Then, we have
    \begin{align}\label{Eq:lem2}
        P\left( \|\widehat{\theta}_S-\theta_S^*\|_2 \leq \frac{5}{\Cmin} \sqrt{d}\lambda_n \right) \geq 1- c_1 p'^{-2},
    \end{align}
    for some constant $c_1>0$.
\end{lemma}
\noindent The last lemma controls the Taylor series remainder term $R^n$:
\begin{lemma}\label{Lem:3}
    If $\lambda_n d \leq \frac{\Cmin^2}{400 \Dmax \DConstThree \log p'}\frac{\alpha}{2-\alpha}$, and $\|W^n\|_{\infty} \leq \frac{\lambda_n}{4}$, then we have
    \begin{align}\label{Eq:lem3}
        P\left( \frac{\|R^n\|_{\infty}}{\lambda_n} \leq \frac{\alpha}{4(2-\alpha)} \right) \geq 1- c_1 p'^{-2},
    \end{align}
    for some constant $c_1>0$.
\end{lemma}

The proof then follows from Lemmas~\ref{Lem:1}, \ref{Lem:2} and \ref{Lem:3} in a straightforward fashion, following \cite{RWL10}. Consider the choice of regularization parameter $\lambda_n = \frac{8(2-\alpha)}{\alpha}\sqrt{\DConstOne\ConstGLM}\sqrt{\frac{\log p}{n}}$. For a sample size greater $n \ge \max\{\frac{4}{\DConstOne \DConstTwo^2 \ConstGLM},\frac{8\ConstHess^2}{\ConstGLM^2}\} \allowbreak \log p$, the conditions of Lemma~\ref{Lem:1} are satisfied, so that we may conclude that  $\|W^n\|_{\infty} \leq \frac{\alpha}{1-\alpha}\frac{\lambda_n}{4} \leq \frac{\lambda_n}{4}$ with high probability. Moreover, with a sufficiently large sample size such that $n \geq L' \left(\frac{2-\alpha}{\alpha}\right)^4 d^2 \DConstOne \DConstThree^2 \log p (\log p')^2$ for some constant $L' > 0$ depending only on $\Cmin$, $\Dmax$, $\ConstGLM$ and $\alpha$, it can be shown that the remaining condition of Lemma \ref{Lem:3} (and hence the milder condition in Lemma \ref{Lem:2}) in turn is satisfied. Therefore, the resulting statements \eqref{Eq:lem2} and \eqref{Eq:lem3} of Lemmas~\ref{Lem:2} and \ref{Lem:3} hold with high probability.

\myparagraph{Strict dual feasibility} Following \cite{RWL10}, we obtain
\begin{align*}
    \|\widehat{Z}_{S^c}\|_{\infty} \ &\leq \ \mNorm{Q^*_{S^cS}(Q^*_{SS})^{-1}}_{\infty} \Big[\frac{\|W_S^n\|_{\infty}}{\lambda_n} + \frac{\|R_S^n\|_{\infty}}{\lambda_n} + 1\Big]
         + \frac{\|W_{S^c}^n\|_{\infty}}{\lambda_n} + \frac{\|R_{S^c}^n\|_{\infty}}{\lambda_n}\\
        &\leq (1-\alpha) +(2-\alpha)\Big[\frac{\|W^n\|_{\infty}}{\lambda_n} + \frac{\|R^n\|_{\infty}}{\lambda_n}\Big]
        \leq (1-\alpha) + \frac{\alpha}{4} + \frac{\alpha}{4} = 1-\frac{\alpha}{2} < 1.
\end{align*}

\myparagraph{Correct sign recovery} To guarantee that the support of $\widehat{\theta}$ is not strictly within the true support $S$, it suffices to show that $\|\widehat{\theta}_S-\theta_S^*\|_{\infty} \leq \frac{\theta_{\min}^*}{2}$. From Lemma~\ref{Lem:2}, we have $\|\widehat{\theta}_S-\theta_S^*\|_{\infty} \leq \|\widehat{\theta}_S-\theta_S^*\|_2 \leq \frac{5}{\Cmin} \sqrt{d}\lambda_n \leq \frac{\theta_{\min}^*}{2}$ as long as $\theta_{\min}^* \geq \frac{10}{\Cmin}\sqrt{d}\lambda_n$. This completes the proof.

\subsection{Proof of Lemma~\ref{Lem:1}}
For a fixed $t \in V\backslash r$, we define $V^{(i)}_t$ for notational convenience so that
\begin{align*}
    W^{n}_{t}  = \frac{1}{n}\sum_{i=1}^n \xs \xt - \xt\dParti(\theta^*_r + \inner{\theta^*_{\backslash r},\xsc}) = \frac{1}{n}\sum_{i=1}^n V^{(i)}_t
\end{align*}
Consider the upper bound on the moment generating function of $V^{(i)}_t$, conditioned on $\xsc$,
\begin{align*}
    &\ \Expect[\exp(a V_t)| \xsc] = \int_{X_r} \exp \bigg\{a \big[X_r \xt - \xt D'\big(\theta^*_r + \inner{\theta^*_{\backslash r},\xsc}\big)\big] \nonumber\\
    &\  + \Big( X_r \big(\theta^*_r + \inner{\theta^*_{\backslash r},\xsc}\big) + C(X_r) -D\big(\theta^*_r + \inner{\theta^*_{\backslash r},\xsc}\big)\Big)\bigg\}\nonumber\\
    = &\ \exp\Big\{ D\big(\theta^*_r + \inner{\theta^*_{\backslash r},\xsc}+a\xt\big) - D\big(\theta^*_r + \inner{\theta^*_{\backslash r},\xsc}\big) - a \xt D'\big(\theta^*_r + \inner{\theta^*_{\backslash r},\xsc}\big)\Big\}\\
    = &\ \exp\Big\{ \frac{a^2}{2}\big({\xt}\big)^2 D''\big(\theta^*_r + \inner{\theta^*_{\backslash r},\xsc}+\nu_i a\xt\big)\Big\} \quad \textrm{for some $\nu_i \in [0,1]$}
\end{align*}
where the last equality holds by the second-order Taylor series expansion. Consequently, we have
\begin{align*}
    &\ \frac{1}{n}\sum_{i=1}^n \log \Expect[\exp(a V^{(i)}_t)| \xsc] \leq \frac{1}{n}\sum_{i=1}^n \frac{a}{2}\big(\xt\big)^2 D''\big(\theta^*_r + \inner{\theta^*_{\backslash r},\xsc}+\nu_i a\xt\big).
\end{align*}
First, we define the event: $\xi_1 := \Big\{\max_{i,s} |\xs| \leq 4 \log p' \Big\}$. Then, by Proposition~\ref{PropSubGaussSingleVariableGLM}, we obtain $P[\xi_1^c] \leq c_1 \, np p'^{-4} \leq c_1 \, p'^{-2}$. Provided that $a \le \DConstTwo$, we can use Condition \ref{AN_GLM_Cond_Smooth} to control the second-order derivative of log-partition function and we obtain
\begin{align*}
    &\ \frac{1}{n}\sum_{i=1}^n \log \Expect[\exp(aV^{(i)}_t)| \xsc] \leq \frac{\DConstOne a^2}{2}\frac{1}{n}\sum_{i=1}^n \big(\xt\big)^2 \quad \textrm{for $a\leq \DConstTwo$}
\end{align*}
with probability at least $1- c_1 p'^{-2}$.
Now, for each index $t$, the variables $\frac{1}{n}\sum_{i=1}^n \big(\xt\big)^2$ satisfy the tail bound in Proposition~\ref{PropSubGaussVariableGLM}. Let us define the event ${\displaystyle\xi_2 := \Big\{\max_{t \in V}} \ \frac{1}{n}\sum_{i=1}^n\big(\xt\big)^2 \allowbreak \leq \ConstGLM \Big\}$ for some constant $\ConstGLM \leq \min\{2\ConstVar/3,2\ConstHess+\ConstVar\}$. Then, we can establish the upper bound of probability $P[\xi_2^c]$ by a union bound,
\begin{align*}
    P[\xi_2^c] \leq \exp (-\frac{\ConstGLM^2}{4 \ConstHess^2} n + \log p) \leq \exp (-c_2 n)
\end{align*}
as long as $n \geq \frac{8\ConstHess^2}{\ConstGLM^2}\log p$. Therefore, conditioned on $\xi_1,\xi_2$, the moment generating function is bounded as follows:
\begin{align*}
    &\ \frac{1}{n}\sum_{i=1}^n \log \Expect[\exp(aV^{(i)}_t)| \xsc,\xi_1,\xi_2] \leq \frac{\DConstOne \ConstGLM \, a^2}{2} \quad \textrm{for $a\leq \DConstTwo$}.
\end{align*}
The standard Chernoff bound technique implies that for any $\delta > 0$,
\begin{align*}
    &\ P\Big[ \big| \frac{1}{n}\sum_{i=1}^n V^{(i)}_t\big | > \delta \ | \ \xi_1,\xi_2 \Big] \leq 2\exp \left(n\Big(\frac{\DConstOne \ConstGLM \, a^2}{2} - a\delta\Big)\right) \quad \textrm{for $a\leq \DConstTwo$}.
\end{align*}
Setting $a = \frac{\delta}{\DConstOne \ConstGLM}$ yields
\begin{align*}
    &\ P\Big[ \big| \frac{1}{n}\sum_{i=1}^n  V^{(i)}_t \big| > \delta \ | \   \xi_1,\xi_2  \Big] \leq 2\exp \left(-\frac{n\delta^2}{2 \DConstOne \ConstGLM}\right) \quad \textrm{for $\delta \leq \DConstOne \DConstTwo \ConstGLM$}.
\end{align*}
Suppose that $\frac{\alpha}{2-\alpha}\frac{\lambda_n}{4} \leq \DConstOne \DConstTwo \ConstGLM$ for large enough $n$; thus setting $\delta = \frac{\alpha}{2-\alpha}\frac{\lambda_n}{4}$:
\begin{align*}
    &\ P\Big[ \big| \frac{1}{n}\sum_{i=1}^n V^{(i)}_t \big| > \frac{\alpha}{2-\alpha}\frac{\lambda_n}{4} \ | \   \xi_1,\xi_2  \Big] \leq 2\exp \left(-\frac{\alpha^2}{(2-\alpha)^2}\frac{n\lambda_n^2}{32\DConstOne \ConstGLM}\right),
\end{align*}
and by a union bound, we obtain
\begin{align*}
    &\ P\Big[ \|W^n\|_{\infty} > \frac{\alpha}{2-\alpha}\frac{\lambda_n}{4} \ | \  \xi_1,\xi_2 \Big] \leq 2\exp \left(-\frac{\alpha^2}{(2-\alpha)^2}\frac{n\lambda_n^2}{32\DConstOne \ConstGLM} + \log p\right).
\end{align*}
Finally, provided that $\lambda_n \geq \frac{8(2-\alpha)}{\alpha}\sqrt{\DConstOne\ConstGLM}\sqrt{\frac{\log p}{n}}$, we obtain
\begin{align*}
    &\ P\Big[ \|W^n\|_{\infty} > \frac{\alpha}{2-\alpha}\frac{\lambda_n}{4} \Big] \leq c_1 p'^{-2} + \exp(-c_2 n) + \exp(-c_3 n),
\end{align*}
where we use the fact that the probability of occurring event $\mathcal{A}$ is upper bounded by $P(\mathcal{A}) \leq P( \xi_1^c) + P(\xi_2^c) + P(\mathcal{A}|\xi_1,\xi_2 )$.

\subsection{Proof of Lemma~\ref{Lem:2}}
In order to establish the error bound $\|\widehat{\theta}_S-\theta_S^*\|_2 \leq B$ for some radius $B$, several works (e.g. \cite{negahban2010,RWL10}) proved that it suffices to show $F(u_S)>0$ for all $u_S := \theta_S-\theta_S^*$ s.t. $\|u_S\|_2 = B$ where
\begin{align*}
    F(u_S) := \ell(\theta_S^* + u_S;X^{1:n})-\ell(\theta_S^*;X^{1:n}) + \lambda_n(\|\theta_S^* + u_S\|_1 - \|\theta_S^*\|_1).
\end{align*}
Note that $F(0) = 0$, and for $\hat{u}_S := \widehat{\theta}_S-\theta_S^*$, $F(\hat{u}_S) \leq 0$. From now on, we show that $F(u_S)$ is strictly positive on the boundary of the ball with radius $B = M\lambda_n\sqrt{d}$ where $M>0$ is a parameter that we will choose later in this proof. Some algebra yields
\begin{align}\label{Eq:lem2tmp}
    F(u_S) \geq (\lambda_n\sqrt{d})^2\Big\{ -\frac{1}{4}M + q^* M^2 - M\Big\}
\end{align}
where $q^*$ is the minimum eigenvalue of $\nabla^2\ell(\theta_S^*+vu_S;X^{1:n})$ for some $v \in [0,1]$. Moreover,
\begin{align*}
    q^* &:= \Lambda_{\min}\big(\nabla^2\ell(\theta_S^*+vu_S)\big)\\
        &\geq \min_{v \in [0,1]} \Lambda_{\min}\big(\nabla^2\ell(\theta_S^*+vu_S)\big)\\
%                   &= \min_{v \in [0,1]} \Lambda_{\min}\big[\frac{1}{n}\sum_{i=1}^n \ddParti(\inner{\theta_S^*+vu_S,\xsc})\xt\big(\xt\big)^T\big]\\
        &\geq \Lambda_{\min}\big[\frac{1}{n}\sum_{i=1}^n \ddParti(\theta^*_r + \inner{\theta_S^*,\xS})\xS\big(\xS\big)^T\big]\\
        & \quad \quad - \max_{v \in [0,1]}\mNorm{\frac{1}{n}\sum_{i=1}^n D'''\big( \theta^*_r
                + \inner{\theta_S^*,\xS}\big)\big(u_S^T\xS\big)\xS\big(\xS\big)^T}_2\\
        &\geq \Cmin -  \max_{v \in [0,1]} \max_{y} \frac{1}{n}\sum_{i=1}^n | D'''\big( \theta^*_r
                + \inner{\theta_S^*,\xS}\big) | \ |\inner{u_S,\xS}| \ \big(\inner{\xS,y}\big)^2
        %&\geq C_{\min} - \Dmax M\lambda_nd \big(\max_{i}|\inner{u_S,\xS}|\big) \max_{v \in [0,1]}\big|\dddParti(\inner{\theta_S^*+vu_S,\xsc})\big|
\end{align*}
where $y \in \reals^d$ s.t $\|y\|_2 =1$.
%    Now define the events
%    \begin{align*}%\label{Eq:event4}
%        \xi_3 &:= \Big\{\max_{i=1,...,n} |\inner{\theta_S^*+vu_S,\xS}| \leq \|\theta^*\|_2 \ \max\{4\ConstMu+2\ConstHessP,9\log p'\} \Big\},\\
%        \xi_1 &:= \Big\{\max_{i,s} |\xs| \leq 4 \log p' \Big\}.
%    \end{align*}
Similarly as in the previous proof, we consider the event $\xi_1$ with probability at least $1-c_1 p'^{-2}$. Then, since all the elements in vector $\xS$ is smaller than $4 \log p'$, $|\inner{u_S,\xS}| \leq 4 \log p' \sqrt{d} \|u_S\|_2 = 4 \log p' M\lambda_n d$ for all $i$. At the same time, by Condition \ref{AN_GLM_Cond_Smooth}, $| \dddParti((\theta^*_r +v u_S) + \inner{\theta_S^*+vu_S,\xS}) | \leq \DConstThree$. Note that $\theta_S^*+vu_S$ is a convex combination of $\theta^*_S$ and $\widehat{\theta}_S$, and as a result, we can directly apply the Condition \ref{AN_GLM_Cond_Smooth}. Hence, conditioned on $\xi_1$, we have
\begin{align*}
    q^* &\geq \Cmin - 4\Dmax M\lambda_n d \DConstThree \log p',
\end{align*}
As a result, assuming that $\lambda_n \leq \frac{\Cmin}{8 \Dmax M d \DConstThree \log p'}$, $q^* \geq \frac{\Cmin}{2}$. Finally, from \eqref{Eq:lem2tmp}, we obtain
\begin{align*}
    F(u_S) \geq (\lambda_n\sqrt{d})^2\Big\{ -\frac{1}{4}M + \frac{\Cmin}{2}M^2 - M\Big\},
\end{align*}
which is strictly positive for $M=\frac{5}{\Cmin}$. Therefore, if $\lambda_n d \leq \frac{\Cmin}{8 \Dmax M \DConstThree \log p'} \leq \frac{\Cmin^2}{40 \Dmax \DConstThree \log p'}$, then
\begin{align*}
    \|\widehat{\theta}_S-\theta_S^*\|_2 \leq \frac{5}{\Cmin} \sqrt{d}\lambda_n,
\end{align*}
which completes the proof.

\subsection{Proof of Lemma~\ref{Lem:3}}
In the proof, we are going to show that $\|R^n\|_{\infty} \leq  4 \DConstThree \log p' \Dmax\|\widehat{\theta}_S-\theta_S^*\|_2^2$. Then, since the conditions of Lemma \ref{Lem:3} are stronger than those of Lemma \ref{Lem:2}, from the result of Lemma \ref{Lem:2}, we can conclude that
\begin{align*}
    \|R^n\|_{\infty} \leq \frac{100 \DConstThree \Dmax \log p'}{\Cmin^2}\lambda_n^2 d,
\end{align*}
as claimed in Lemma \ref{Lem:3}.

From the definition of $R^n$, for a fixed $t \in V\backslash r$, $R_t^n$ can be written as
\begin{align*}
    \frac{1}{n}\sum_{i=1}^n \Big[D''\big(\theta^*_r + \inner{\theta^*_{\backslash r},\xsc}\big) - D''\big(\bar{\theta}_s + \inner{\bar{\theta}^{(t)},\xsc}\big)\Big]\Big[\xsc\big(\xsc\big)^T\Big]_t^T[\widehat{\theta}_{\backslash r}-\theta^*_{\backslash r}]
\end{align*}
where $\bar{\theta}^{(t)}_{\backslash r}$ is some point in the line between $\widehat{\theta}_{\backslash r}$ and $\theta^*_{\backslash r}$, i.e., $\bar{\theta}^{(t)}_{\backslash r} = v_t\widehat{\theta}_{\backslash r} + (1-v_t)\theta^*_{\backslash r}$ for $v_t \in [0,1]$. By another application of the mean value theorem, we have
\begin{align*}
    R_t^n = -\frac{1}{n}\sum_{i=1}^n \left\{D'''\big(\bar{\bar{\theta}}_s + \inner{\bar{\bar{\theta}}^{(t)}_{\backslash r},\xsc}\big)\xt\right\}\left\{v_t[\widehat{\theta}_{\backslash r}-\theta^*_{\backslash r}]^T\xsc\big(\xsc\big)^T[\widehat{\theta}_{\backslash r}-\theta^*_{\backslash r}]\right\}
\end{align*}
for a some point $\bar{\bar{\theta}}^{(t)}_{\backslash r}$ between $\bar{\theta}^{(t)}_{\backslash r}$ and $\theta^*_{\backslash r}$. Similarly in the previous proofs, conditioned on the event $\xi_1$, we obtain
\begin{align*}
    |R_t^n| \leq \frac{4 \DConstThree \log p'}{n} \sum_{i=1}^n \left\{v_t[\widehat{\theta}_{\backslash r}-\theta^*_{\backslash r}]^T\xsc\big(\xsc\big)^T[\widehat{\theta}_{\backslash r}-\theta^*_{\backslash r}]\right\}.
\end{align*}
Performing some algebra yields
\begin{align*}
    |R_t^n| \leq 4 \DConstThree \Dmax  \log p' \, \|\widehat{\theta}_S-\theta_S^*\|_2^2, \quad \textrm{for all $t \in V\backslash r$}
\end{align*}
with probability at least $1 - c_1 p'^{-2}$, which completes the proof.

\section{Optimization Problems for Poisson and Exponential Graphical
Model Neighborhood Selection}\label{appendix:optimization}

We propose to fit our family of graphical models by performing
neighborhood selection, or maximizing the $\ell_{1}$-penalized
log-likleihood for each node conditional on all other nodes.  For
several exponential families, however, further restrictions on the
parameter space are needed to ensure a proper Markov Random Field.
When performing neighborhood selection, these can be imposed by adding
constraints to the penalized generalized linear models.  We illustrate
this by providing the optimization problems solved by our Poisson
graphical model and exponential graphical model M-estimator that are
used in Section 4.

Following from Section 3, the neighborhood selection problem for our
family of models maximizes the likelihood of a node, $X_{r}$,
conditional on all other nodes, $X_{V \backslash r}$.  This
conditional likelihood is regularized with an $\ell_{1}$ penalty to
induce sparsity in the edge weights, $\theta(r)$ a $p-1$ dimensional
vector, and constrained to
enforce restrictions, $\theta(r) \in \mathcal{C}$, needed to yield a
proper MRF:
\begin{align*}
\maximize_{\theta(r)} \ \ \frac{1}{n} \sum_{i=1}^{n} \ell \left(
X_{i,r} | X_{i, V \backslash r} ; \theta(r) \right) - \lambda_{n} \|
\theta(r) \|_{1}  \ \ \textrm{subject to} \ \ \theta(r)
\in \mathcal{C},
\end{align*}
where $\ell \left( X_{i,r} | X_{i, V \backslash r} ; \theta(r)
\right)$ is the conditional log-likelihood for the exponential
family.  For the Poisson graphical model, the edge weights are
constrained to be non-positive.  This yields the following
optimization problem:
\begin{align*}
\maximize_{\theta(r)} & \ \ \frac{1}{n} \sum_{i=1}^{n} \left[
X_{r,i} X_{V \backslash r,i}^{T} \theta(r) - \mathrm{exp} \left( X_{V
\backslash r,i}^{T} \theta(r) \right) \right]
- \lambda_{n} \| \theta(r) \|_{1} \\
\textrm{subject to } & \ \ \theta(r) \leq 0.
\end{align*}
Similarly, the edge weights of the exponential graphical are
restricted to be non-negative yielding:
\begin{align*}
\maximize_{\theta(r)} & \ \ \frac{1}{n} \sum_{i=1}^{n} \left[ -
X_{r,i} X_{V \backslash r,i}^{T} \theta(r) + \mathrm{log} \left( X_{V
\backslash r,i}^{T} \theta(r) \right) \right]
- \lambda_{n} \| \theta(r) \|_{1} \\
\textrm{subject to } & \ \ \theta(r) \geq 0.
\end{align*}
Note that we neglect the intercept term, assuming this to be zero as
is common in other proposed neighborhood selection methods
\citep{Meinshausen06,RWL10,JRVS11}.
Both of the neighborhood selection problems are concave problems with
a smooth log-likelihood and linear constraints.
While there are a plethora of optimization routines available to solve
such problems, we have employed a projected gradient descent scheme
which is guaranteed to converge to a global optimum
\citep{daubechies2008accelerated,beck_teb_2010}.

% Note: in this sample, the section number is hard-coded in. Following
% proper LaTeX conventions, it should properly be coded as a reference:

%In this appendix we prove the following theorem from
%Section~\ref{sec:textree-generalization}:

\vskip 0.2in
\bibliography{glmgm,network,sml}

\end{document}